\documentclass[10pt]{article}
\usepackage{amsmath}
\usepackage{latexsym}
\usepackage{amssymb}
\usepackage{amsfonts}
\usepackage{amsthm}
\title{GENERATING COMPATIBILITY CONDITIONS   \\ IN MATHEMATICAL PHYSICS  }
\author{J.-F. Pommaret \\ CERMICS, Ecole des Ponts ParisTech, 77455 France \\
 jean-francois.pommaret@wanadoo.fr \\
 http://cermics.enpc.fr/$\sim$pommaret/home.html }
\date{  }
\begin{document}
\maketitle

\vspace{2cm}

\noindent
{\bf ABSTRACT}  \\

The search for {\it generating compatibility conditions} (CC) for a given operator is a very recent problem met in General Relativity in order to study the {\it Killing} operator for various standard useful metrics (Minkowski, Schwarschild and Kerr). In this paper, we prove that the link existing between the lack of {\it formal exactness} of an operator sequence on the jet level, the lack of {\it formal exactness} of its corresponding symbol sequence and the lack of {\it formal integrability} (FI) of the initial operator is of a purely homological nature as it is based on the {\it long exact connecting sequence} provided by the so-called {\it snake lemma}. It is therefore quite difficult to grasp it in general and even more difficult to use it on explicit examples. It does not seem that any one of the results presented in this paper is known as most of the other authors who studied the above problem of computing the total number of generating CC are confusing this number with a kind of {\it differential transcendence degree}, also called {\it degree of generality} by A. Einstein in his 1930 letters to E. Cartan. The motivating examples that we provide are among the rare ones known in the literature and could be used as testing examples for future applications of computer algebra.   \\

\vspace{3cm}


\newpage   
\noindent
{\bf 1) INTRODUCTION}  \\

If $X$ is a manifold of dimension $n$ with local coordinates $(x^1,...,x^n)$, let us introduce the tangent bundle $T=T(X)$ and the cotangent bundle $T^*=T^*(X)$, the $q$-symmetric tensor bundle $S_qT^*$ and the bundle ${\wedge}^rT^*$ of $r$-forms. In General Relativity, there may be different solutions of Einstein equations in vacuum like the Minkowski, the Schwarzschild and the Kerr metrics for example. For fixing the notations and with more details, if $\omega \in S_2T^*$ is a nondegenerate metric, that is $det(\omega)\neq 0$, and if $j_q$ denotes all the derivatives of an object up to order $q$, we may construct the Christoffel symbols $\gamma$ through the Levi-Civita isomorphism $(\omega, \gamma) \simeq  j_1(\omega)$ and, using the language of jet bundles, $(\omega,\gamma)$ is a section of $J_1(S_2T^*)$ that will be simply written $(\omega, \gamma)\in J_1(S_2T^*)$. Then we can introduce the well-known Riemann tensor $\rho=({\rho}^k_{l,ij})\in {\wedge}^2T^*\otimes T^*\otimes T$ with ${\rho}_{kl,ij}+{\rho}_{lk,ij}=0$ after lowering the upper index by means of $\omega$ {\it and} $\delta \rho=0$ where $\delta : {\wedge}^2T^*\otimes T^* \otimes T \rightarrow {\wedge}^3T^*\otimes T$ is the Spencer $\delta$-map. Introducing the Ricci tensor ${\rho}_{ij}= {\rho}^r_{i,rj}={\rho}_{ji}$ or the Einstein tensor ${\epsilon}_{ij}={\rho}_{ij}-\frac{1}{2}{\omega}_{ij}{\omega}^{rs}{\rho}_{rs}$, the $10$ non-linear Einstein equations are described by ${\epsilon}_{ij}=0$ or, equivalently, by ${\rho}_{ij}=0$ when $n=4$. \\

Now, if ${\cal{E}}$ is a fibered manifold over $X$ with fiber dimension $m$ and local coordinates $(x^i,y^k)$ with $ i=1,...,n$ and $k=1,... ,m$, we may introduce the tangent bundle $T({\cal{E}})$ over ${\cal{E}}$ with local coordinates $(x,y,u,v)$ and the {\it vertical bundle} $V({\cal{E}})$ with local coordinates $(x,y,u=0,v)=(x,y,v)$ which are both vector bundles over ${\cal{E}}$. We shal dnote by the capital letters $\Omega \in S_2T^*, \Gamma \in S_2T^*\otimes T, R\in {\wedge}^2T^*\otimes T^*\otimes T$, the respective {\it linearizations} of $\omega, \gamma, \rho$ which are sections of the respective vertical bundles. Introducing the {\it Lie } derivative ${\cal{L}}$ of geometric objects, it is therefore possible to introduce the corresponding first order {\it Killing} operator ${\cal{D}}:T \rightarrow \Omega \in S_2T^*: \xi \rightarrow {\cal{L}}(\xi)\omega$, the first order {\it Christoffel} operator $\Omega \rightarrow \Gamma\in S_2T^*\otimes T$ in such a way that $Christoffel \circ  Killing :  T \rightarrow \Gamma \in S_2T^*\otimes T: \xi \rightarrow {\cal{L}}(\xi)\gamma$ and the second order {\it Riemann} operator $ \Omega \rightarrow R\in {\wedge}^2T^*\otimes T^*\otimes T$ in such a way that $Riemann  \circ Killing : T \rightarrow R \in {\wedge}^2T^*\otimes T^*\otimes T:\xi \rightarrow {\cal{L}}(\xi)\rho$ both with its {\it contraction} $\Omega \rightarrow S_2T^*$ called {\it Ricci} operator. For example, it is known that $ 2\,{\omega}_{rk} {\gamma}^k_{ij}= ({\partial}_i{\omega }_{rj} +{\partial} _j{\omega}_{ir} - {\partial}_r{\omega}_{ij})$ that we shall write simply, using formal notations, 
$ 2 \, \omega \gamma= ( \partial \omega + \partial \omega - \partial \omega)$ and thus $ 2 \omega \Gamma +2 \gamma  \Omega =(d\Omega + d\Omega  - d \Omega) $.                                                          
We have proved in ([16],[18],[19],[20]) that the so-called gravitational waves equations are nothing else than $ad(Ricci)$ by introducing the {\it formal adjoint} operator. It is important to notice that the {\it Einstein } operator $\Omega \rightarrow E_{ij}= R_{ij} - \frac{1}{2} {\omega}_{ij}{\omega}^{rs}R_{rs}$ is self-adjoint with $6$ terms though the {\it Ricci} operator is {\it not} with only $4$ terms. Recently, many physicists (See [1],[2],[8],[9],[24]) have tried to construct the {\it compatibility conditions} (CC) of the {\it Killing } operator for various types of background metrics, in particular the three ones already quoted, namely an operator ${\cal{D}}_1:S_2T^* \rightarrow F_1$ such that ${\cal{D}}_1\Omega=0$ generates the CC of ${\cal{D}}\xi=\Omega$. We have proved in the above references the following crucial results:  \\

\noindent
$\bullet$ These CC may contain a certain number of second and third order CC. It is therefore crucial in actual practice to select the successive generating CC of order $1,2,3,...$ till we stop because of noetherian arguments ([14]).   \\

\noindent
$\bullet$ These CC only depend on the Lie algebra structure (dimension of the solution space and commutation relations) of the corresponding Killing operator, which, even though it is finite dimensional with dimension $\leq n(n+1)/2$ that is $10$ obtained for the Minkowski metric, may have dimension $4$ for the Schwarzschild metric and dimension $2$ for the Kerr metric.   \\

\noindent
$\bullet$ The only two canonical sequences that can be constructed from an operator or a system, namely the Janet and Spencer sequences, are structurally quite different. Indeed, the {\it Janet bundles} $F_0, ... F_n$ appearing in the Janet sequence are concerned wit geometric objects like $\omega, \gamma, \rho$, while the {\it Spencer bundles} $C_0, ... ,C_n$ are far from being related with geometric objects, the simplest example being $C_0=R_q\subseteq J_q(E)$. In the case of Lie equations considered, the central concept is not the {\it system} but rather the {\it group} as it can be seen at once from the construction of the {\it Vessiot structure equations} ([12],[13],[17],[19]).   \\

The authors wo have studied these questions had in mind that the total number of generating CC could be considered as a kind of  "{\it differential transcendence degree} ", also called "{\it degree of generality} " by A. Einstein in his letters to E. Cartan of $1930$ on absolute parallelism ([3]), the modern definition being that of  the  "{\it differential rank} " ([13],[14],[23]). We must say that Cartan, being unable to explain to Einstein his theory of exterior systems, just copied the work of Janet in his letters to Einstein, published later on as the {\it only} paper he wrote on the PD aproach, but without ever quoting Janet who suffered a lot from this behaviour and had to turn to mechanics.  \\

 Such a result will be obtained in the framework of differential modules as its explanation in the framework of differential systems is much more delicate and technical ([7],[15],[23]). \\
First of all, with our previous assumptions, $D=K[d]$ is a noetherian domain and we can restrict our study to finitely generated differential modules which are therefore finitely presented (See [14] for more details). Let thus $M$ be defined by a finite free presentation giving rise to the long exact sequence:  \\
\[     0 \rightarrow L \rightarrow D^p \stackrel{{\cal{D}}}{\longrightarrow } D^m \stackrel{p}{\longrightarrow} M \rightarrow 0  \]
where the differential operator ${\cal{D}}$ is acting on the right by composition with action law $(P,{\cal{D}}) \rightarrow P \circ {\cal{D}}, \forall P\in D$, $p$ is the canonical residual projection and $L=ker({\cal{D}})\subset D^p$. The image $im({\cal{D}})\subset D^m$ is called the {\it differential module of equations} and is thus finitely generated because $D$ is a noetherian differential domain. \\

\noindent
{\bf DEFINITION 1.1}: The {\it differential rank} $rk_D(M)$ over $D$ of a differential module $M$ is the differential rank over $D$ of the maximum free differential submodule $F$ of $M$ and we have the short exact sequence $0 \rightarrow F \rightarrow M \rightarrow T \rightarrow 0$ where $T=M/F$ is a torsion module over $D$. In particular, if $F\simeq D^r$, then $rk_D(M)=r$.    \\

The following quite useful proposition proves the {\it additivity property } of the differential rank and is used in the next corollary ([11],[13],[22]):  \\

\noindent
{\bf PROPOSITION 1.2}: If we have a short exact sequence $0 \rightarrow M' \stackrel{f}{\rightarrow } M \stackrel{g}{\rightarrow } M" \rightarrow 0$ of differential modules, then $rk_D(M)=rk_D(M') + rk_D(M")$.  \\

\noindent
{\bf COROLLARY 1.3}: ({\it Euler-Poincar\'{e} characteristic}) For {\it any} finite free differential resolution of a differential module $M$, then $rk_D(M)$ is equal to the alternate sum ${\chi}_D(M)$ of differential ranks of the free differential modules of the resolution.   \\  \\

We obtain therefore $rk_D(L) -p +m -rk_D(M)=0$ and it follows from noetherian arguments that the differential module $L\subset D^p$ {\it is finitely generated but not free in general} and we may look for a minimum number of generators which may be differentially dependent in general as we shall see in the next examples. It thus remains to provide examples of such computations showing that these two numbers are not related and must therefore be found totally independently in general, apart from the very exceptional situation met when there is only a single generating CC.  \\ 

In actual practice, working in the system framework, starting with a system $R_q\subset J_q(E)$ of order $q$ on $E$ and introducing the {\it canonical projection} $\Phi: J_q(E) \rightarrow  F_0=J_q(E)/R_q$, we shall construct for each $r\geq 0$ a family of FI systems $B_r = im ({\rho}_r(\Phi)) \subseteq J_r(F_0)$ such that $B_{r+1}\subseteq {\rho}_1(B_r)=J_1(B_r)\cup J_{r+1}(F_0)\subset J_1(J_r(F_0))$ projects onto $B_r$, that is $B_{r+1}$ is defined by more generating PD equations than the ones defining $B_r$ both with its prolongations, and start to get equality when $r$ is large enough in the projective limit $B_{\infty} \rightarrow ... \rightarrow B_{r+1} \rightarrow B_r \rightarrow ... \rightarrow B_1 \rightarrow F_0 \rightarrow 0$ . {\it The striking result is that there may be gaps in the procedure}, that is we shall even provide a tricky example where one can have a single generating CC of order $3$, then no {\it new} generating CC of order $4$ {\it and} $5$, but {\it suddenly} a new generating CC of order $6$ ending the procedure. We do not believe that such situations were even known to exist. \\

\newpage
\noindent
{\bf 2) MOTIVATING EXAMPLES}  \\

We provide below three examples, pointing out that it is quite difficult to exhibit such examples.\\

\noindent
{\bf EXAMPLE 2.1}: With $n=3, m=dim(E)=2, dim (F_0)=5 $ and $K=\mathbb{Q}$ while keeping an upper index for any unknown, let us consider the following system $R_1\subset J_1(E)$ with $dim(R_1)=3$ because $par_1=\{{\xi}^1,{\xi}^2, {\xi}^2_1\}$ and corresponding Janet tabular:  \\
\[ \left\{   \begin{array}{rcl}
{\Phi}^5 & \equiv & {\xi}^2_3  =0  \\
{\Phi}^4 & \equiv & {\xi}^1_3 =0  \\
{\Phi}^3 & \equiv & {\xi}^2_2=0  \\
{\Phi}^2 & \equiv & {\xi}^1_2 + {\xi}^2_1=0    \\
{\Phi}^1 & \equiv & {\xi}^1_1=0  
\end{array}  \right.  \fbox{  $  \begin{array} {ccc}
1 & 2 & 3  \\
1 & 2 & 3   \\
1 & 2 & \bullet  \\
1 & 2 & \bullet  \\
1 & \times & \bullet  
\end{array}  $  }  \] 
It is easy to check that all the second order jets vanish and that the general solution $\{ {\xi}^1= ax^2 + b, {\xi}^2= cx^1 + d \mid a+c=0 \}$ depends on $3$ arbitrary constants. As the non-multiplicative variable written with the sign $\times$ cannot be used, the symbol $g_1$ is not involutive because it is finite type with $g_2=0$. This system is trivially FI because it is made by homogeneous PD equations. We have the following commutative diagrams:  

 \[  \begin{array}{rccccccccccl}
 &&&& 0 &&0&&0&  & \\
 &&&& \downarrow && \downarrow && \downarrow  & & \\
  && 0 & \longrightarrow &S_3T^*\otimes E&\longrightarrow&S_2T^*\otimes F_0 &\longrightarrow & h_2 & \longrightarrow  0  \\
   &&\downarrow && \downarrow & & \downarrow & & \downarrow & &    & \\
0& \longrightarrow &R_3&\longrightarrow &J_3(E)& \longrightarrow &J_2(F_0) &\longrightarrow & Q_2 &   \longrightarrow 0 \\
    & & \downarrow && \hspace{5mm}\downarrow {\pi}^3_2 & &\hspace{5mm} \downarrow {\pi}^2_1 & & \downarrow  &  \\
   0 & \longrightarrow & R_2& \longrightarrow & J_2(E) & \longrightarrow &J_1(F_0)& \longrightarrow &Q_1& \longrightarrow  0 \\
   && \downarrow && \downarrow & & \downarrow & & \downarrow &      \\
   &&0&& 0 && 0 && 0  &  
   \end{array}     \]
 \[  \begin{array}{rccccccccccl}
 &&&& 0 &&0&& 0 &  & \\
 & &&& \downarrow && \downarrow &&\downarrow &  & \\
  && 0 & \longrightarrow &20 &\longrightarrow& 30 &\longrightarrow & 10 &\longrightarrow  0  \\
   &&\downarrow && \downarrow & & \downarrow & & \downarrow & &    & \\
0& \longrightarrow &3&\longrightarrow & 40 & \longrightarrow & 50 &\longrightarrow & 13 &   \longrightarrow 0 \\
    & & \downarrow && \hspace{5mm}\downarrow {\pi}^3_2 & &\hspace{5mm} \downarrow {\pi}^2_1 & & \downarrow  &  \\
   0 & \longrightarrow & 3 & \longrightarrow & 20 & \longrightarrow & 20 & \longrightarrow &  3  & \longrightarrow  0 \\
   &&\downarrow && \downarrow & & \downarrow & & \downarrow &      \\
   && 0 && 0 && 0 && 0  &  
   \end{array}     \]

 \[  \begin{array}{rccccccccccl}
 &&&& 0 &&0&&&  & \\
 &&&& \downarrow && \downarrow && & & \\
  && 0 & \rightarrow &S_3T^*\otimes E&\rightarrow&S_2T^*\otimes F_0 &\rightarrow & h_2 &
   \rightarrow  0  \\
   &&&& \hspace{3mm} \downarrow \delta & & \hspace{3mm} \downarrow \delta & & \downarrow & &    & \\
&  &0&\rightarrow &T^*\otimes S_2T^*\otimes E & \rightarrow &T^*\otimes T^* \otimes F_0 &\rightarrow & T^*\otimes h_1 &   \rightarrow 0 \\
    & & && \hspace{3mm}\downarrow \delta  & &\hspace{3mm} \downarrow \delta& & \downarrow  &  \\
   0 & \rightarrow & {\wedge}^2T^*\otimes g_1& \rightarrow & {\wedge}^2T^*\otimes T^* \otimes E & \rightarrow & {\wedge}^2T^*\otimes F_0  & \rightarrow & 0 &  \\
   && \hspace{3mm} \downarrow \delta && \hspace{3mm} \downarrow \delta & & \downarrow & & &   &   \\
 0   & \rightarrow & {\wedge}^3T^*\otimes E  & = & {\wedge}^3T^*\otimes E & \rightarrow & 0  & &&&   \\
       &   &  \downarrow &  & \downarrow  &  &  &  &  &  &  \\
   &&0&& 0 && && &  &
   \end{array}     \]

 \[  \begin{array}{rccccccccccl}
 &&&& 0 &&0&&&  & \\
 &&&& \downarrow && \downarrow && & & \\
  && 0 & \rightarrow &20&\rightarrow&30 &\rightarrow & 10 &
   \rightarrow  0  \\
   &&&& \hspace{3mm} \downarrow \delta & & \hspace{3mm} \downarrow \delta & & \downarrow & &    & \\
&  &0&\rightarrow &36 & \rightarrow &45 &\rightarrow & 9 &   \rightarrow 0 \\
    & & && \hspace{3mm}\downarrow \delta  & &\hspace{3mm} \downarrow \delta& & \downarrow  &  \\
   0 & \rightarrow & 3 & \rightarrow & 18 & \rightarrow & 15 & \rightarrow & 0 &  \\
   && \hspace{3mm} \downarrow \delta && \hspace{3mm} \downarrow \delta & & \downarrow & & &   &   \\
 0   & \rightarrow &  2  & = &  2  & \rightarrow & 0  & &&&   \\
       &   &  \downarrow &  & \downarrow  &  &  &  &  &  &  \\
   &&0&& 0 && && &  &
   \end{array}     \]

The next result points out the importance of the Spencer $\delta$-cohomology. Indeed, we shall prove that the last symbol diagram is commutative and exact. In particular, the lower left map $\delta$ is surjective and thus  the upper right induced map $h_2 \rightarrow T^*\otimes Q_1$ is  also surjective while these two maps have isomorphic kernels. \\
For this, we notice that the $3$ components of ${\wedge}^2T^*\otimes g_1$ are $\{ v^2_{1, 12}, v^2_{1,13}, v^2_{1,23}\}$ and the map $\delta$ is described by the two linear equations:  \\
\[  w^1_{123}\equiv v^1_{1,23} +v^1_{2,31} + v^1_{3,12}=v^2_{1,13}=0, \hspace{5mm}  w^2_{123}\equiv v^2_{1,23} +v^2_{2,31} + v^2_{3,12}=v^2_{1,23}=0 \]
that is to say by two linearly independent equations. Accordingly, in the left column we have:  \\
\[    dim(H^2(g_1))=dim(Z^2(g_1))=dim(ker(\delta))=1\]   
An unusual snake-type diagonal chase left to the reader as an exercise proves that the induced map $h_2\rightarrow T^*\otimes Q_1$ is surjective with a kernel isomorphic to $H^2(g_1)$. This is indeed a {\it crucial result} because it also proves that the additional CC has only to do with the the single second order component of the Riemann tensor in dimension $2$, a striking result that could not even be imagined by standard methods. Moreover, we know that if a system $R_q\subset J_q(E)$ is FI, for example when it is homogeneous like in this case,  and its symbol $g_q \subset S_qT^*\otimes E$ is such that $s$ is the smallest integer such that $g_{q+s}$ becomes $2$-acyclic (or involutive), then the generating CC are of order at most $s+1$ ([13],[14],[19]).  \\
 
Collecting the above results, we find the $3$ {\it first order} differentially independent generating CC coming from the Janet tabular and the {\it additional} single {\it second order} generating CC describing the $2$-dimensional {\it Riemann} operator, that is the linearized Riemann tensor in the space $(x^1,x^2)$:  \\

\[  \left\{     \begin{array}{rcl}
{\Psi}^4 & \equiv & d_{22}{\Phi}^1 + d_{11}{\Phi}^3 - d_{12}{\Phi}^2  =0  \\
              &           &                        \\
{\Psi}^3 & \equiv &   d_3{\Phi}^3 - d_2{\Phi}^5=0    \\
{\Psi}^2 & \equiv &   d_3{\Phi}^2 - d_2{\Phi}^4 - d_1{\Phi}^5  =0        \\
{\Psi}^1 & \equiv & d_3{\Phi}^1 - d_1{\Phi}^4  =0
\end{array} \right.  \] 

An elementary computation provides the second order CC:  \\
\[   d_{22}{\Psi}^1 +d_{11}{\Psi}^3 - d_{12}{\Psi}^2 -  d_3{\Psi}^4 =0   \]

The corresponding differential sequence written with differential modules over the ring $D=K[d_1,d_2,d_3]=K[d]$ is:  \\
\[          0 \rightarrow D  \underset{2}{\rightarrow} D^4 \underset{2}{\rightarrow} D^5 \underset{1}{\rightarrow} D^2  \stackrel{p}{\rightarrow}M 
\rightarrow 0   \]
where $p$ is the canonical (residual) projection. We check indeed that $ 1 - 4 + 5 -2 = 0$ but this sequence is quite far from being even strictly exact. Of course, as $R_2$ is involutive, we may set 
$C_r={\wedge}^rT^*\otimes R_2$ and obtain the corresponding canonical second Spencer sequence which is induced by the Spencer operator:\\ 
\[   0  \rightarrow \Theta \stackrel{j_2}{\longrightarrow} C_0 \stackrel{D_1}{\longrightarrow} C_1 \stackrel{D_2}{\longrightarrow}  C_2 \stackrel{D_3}{\longrightarrow}C_3 \rightarrow 0    \]  with dimensions:  \\
\[   0  \rightarrow \Theta \stackrel{j_2}{\longrightarrow} 3 \underset{1}{\stackrel{D_1}{\longrightarrow}} 9 \underset{1}{\stackrel{D_2}{\longrightarrow}}  9 \underset{1}{ \stackrel{D_3}{\longrightarrow}} 3 \rightarrow 0    \]
Proceeding inductively as we did for finding the second order CC, we may obtain by combinatorics  the following formally exact sequence:  \\
\[  0 \rightarrow  \Theta \rightarrow 2 \underset 1{\rightarrow} 5 \underset 2{\rightarrow}13 \underset 1{\rightarrow} 19 \underset 1{\rightarrow} 12 \underset 1{\rightarrow} 3 \rightarrow 0  \]
with Euler-Poincar\'{e} characteristic $2 - 5 + 13 - 19 + 12 - 3=0$ but, as before, there is a matrix $260 \times 280$ at least and we doubt about the use of computer algebra, even on such an elementary example. With $F_0=J_1(E)/R_1$, the starting long exact sequence used as a middle row of the first diagram with dimensions:  \\
\[     0 \rightarrow  3  \rightarrow 40 \rightarrow 50 \rightarrow 13 \rightarrow 0  \]
and we have $13=(3 + 3\times 3) +1$, that is three generating first order CC which are differentially independent, plus their $9$ prolongations, plus one second order CC which is nevertheless {\it not} differentially independent. Hence we have a total number of $3+1=4$ generating CC but this number has nothing to do with any differential transcendence degree because ${\Psi}^4$ is differentially algebraic over $\{ {\Psi}^1,{\Psi}^2,{\Psi}^3\}$. \\
We finally compute the corresponding (canonical) Janet sequence by quotient. For this, we must use the trivial second Spencer sequence:  \\
\[   0  \rightarrow E \stackrel{j_2}{\longrightarrow} C_0(E) \stackrel{D_1}{\longrightarrow} C_1(E) \stackrel{D_2}{\longrightarrow}  C_2(E) \stackrel{D_3}{\longrightarrow}C_3(E) \rightarrow 0    \] 
namely:   \\
\[  0 \rightarrow  2 \underset 2{\rightarrow} 20 \underset 1{\rightarrow} 40 \underset 1{\rightarrow}30 \underset 1{\rightarrow} 8  \rightarrow 0  \]
with $2 - 20 + 40 - 30 + 8=0$. The (canonical) Janet sequence is thus:  \\
\[  0 \rightarrow  \Theta \rightarrow E \underset 2{\rightarrow} F_0 \underset 1{\rightarrow}F_1 \underset 1{\rightarrow} F_2 \underset 1{\rightarrow} F_3 \rightarrow 0  \]
with now $F_0=J_2(E)/R_2$ and $F_r=C_r(E) / C_r , \forall r = 0,1,2,3$ and dimensions:   \\
\[  0 \rightarrow  \Theta \rightarrow 2 \underset 2{\rightarrow} 17 \underset 1{\rightarrow}31 \underset 1{\rightarrow} 21 \underset 1{\rightarrow} 5  \rightarrow 0  \]
so that we have again $2 - 17 + 31 - 21 +5=0$ in a coherent way with the fact that $rk_D(M)=0$.     \\  \\

\noindent
{\bf EXAMPLE 2.2}: With $n=3,m=1, q=2$ and $K=\mathbb{Q}(x^1,x^2,x^3)=\mathbb{Q}(x)$, let us consider the following linear inhomogeneous system:   \\
\[        y_{33} - x^2y_1=v,  \,\,\,\,\,  y_{12}=u  \]

\noindent
$\bullet$ {\it Step} 1: The symbol $g_2$ is defined by $v_{33}=0, v_{12}=0$ may not be involutive {\it or} the coordinate system may not be $\delta$-regular. However, changing linearly the local coordinates with $x^1 \rightarrow x^1, x^2\rightarrow x^2+x^1, x^3 \rightarrow x^3$, we obtain the Janet tabular for $g_2$:  \\
\[ \left \{ \begin{array}{lcl}
  v_{33} &  =  & 0  \\
  v_{22} + v_{12} &  =  &  0 \\
  \end{array} \right. \fbox{ $\begin{array}{lll}
  1 & 2 & 3  \\
  1 & 2 & \bullet  
  \end{array} $}   \]
  and thus the Janet tabular for $g_3$:  \\
  \[ \left \{ \begin{array}{lcl}
  v_{333} &  =  & 0  \\
  v_{233}  & = & 0   \\
  v_{223} + v_{123} &  =  &  0 \\
  v_{222} + v_{122}  &  =  &  0  \\
  v_{133}  &  =  &  0  \\
  v_{122} + v_{112}  &  =  & 0
  \end{array} \right. \fbox{ $\begin{array}{ccc}
  1 & 2 & 3  \\
  1  &  2  &  \bullet  \\
  1  &  2  &  \bullet  \\
  1 & 2 & \bullet   \\
  1  &  \bullet  &  \bullet  \\
  1  &  \bullet  &  \bullet  \\
  \end{array} $ }   \] 
 We let the reader check as an exercise that $g_2$ is not $2$-acyclic by counting the dimensions in the long sequence: \\
 \[  0 \rightarrow g_4  \stackrel{\delta}{\rightarrow} T^* \otimes g_3 \stackrel{\delta}{\rightarrow} {\wedge}^2T^*\otimes g_2 \stackrel{\delta}{\rightarrow} {\wedge}^3T^*\otimes T^*\]
and that $g_3$ is involutive, thus $2$-acyclic, with characters $(0,0,4)$. It follows that $dim(g_2)=6-2=4$, $dim(g_3)=0+0+4=4, dim(g_4)=4, ...$. We obtain from the main theorem ${\rho}_r(R^{(1)}_2)= R^{(1)}_{2+r}$. It is easy to chek that $R^{(1)}_2=R_2$ with $dim(R_2)=8$, $dim(R_3)=8+4=12$, $dim(R_4)=11+4=15$, $im(R_5)=56-39=13+4=17$ but things are changing after that. As such a property is intrinsic, coming back to the original system of coordinates, we have after one more prolongation:   \\
 \[ \left \{ \begin{array}{lcl}
 \,\,\,\,\,y_{1233}  & = &  u_{33}   \\
 -y_{1233}+x^2y_{112}+ y_{11}  &  =  &  -v_{12}  \\
 -x^2y_{112}  &  =  &  - x^2u_1 
  \end{array}\right. \]
 and thus $y_{11}=u_{33}- v_{12} - x^2u_1$. We may thus consider the new second order system $
 R'_2=R^{(2)}_2\subset R_2$ with a strict inclusion and $dim(R^{(2)}_2)=7$:    \\
\[        y_{33} - x^2y_1=v,  \,\,\,\,\,  y_{12}=u , \,\,\,\, y_{11}=w, \,\,\,\,  w=u_{33}-v_{12} - x^2 u_1 \]

We may start again with $R'_2$ and study its symbol $g'_2 $ defined by the 3 linear equations with the following Janet tabular obtained after doing the same change of local coordinates as before:  \\
\[  \left \{ \begin{array}{lcl}
  v_{33} &  =  & 0  \\
  v_{22} + v_{12} &  =  &  0 \\
  v_{12} + v_{11} &  =  &  0
  \end{array} \right. \fbox{ $\begin{array}{ccc}
  1 & 2 & 3  \\
  1 & 2 & \bullet  \\
  1  &  \bullet & \bullet   
\end{array} $ }   \]
This symbol is neither $2$-acyclic nor involutive but its prolongation $g'_3$, defined by the $8$ equatios:  \\
\[  \left \{ \begin{array}{lcl}
  v_{333} &  =  & 0  \\
  v_{233}  & =  & 0   \\
  v_{223} - v_{113} &  =  &  0 \\
  v_{222} + v_{111} &  =  &  0  \\
  v_{133} & = & 0  \\
  v_{123} + v_{113} & = & 0  \\
  v_{122} -v_{111} & = & 0  \\
  v_{112} + v_{111} & = & 0
  \end{array} \right. \fbox{ $\begin{array}{ccc}
  1 & 2 & 3  \\
  1 & 2 & \bullet  \\
  1 & 2 & \bullet \\
  1 & 2 & \bullet  \\
  1 & \bullet & \bullet  \\
  1 & \bullet & \bullet  \\
  1 & \bullet & \bullet  \\
  1  &  \bullet & \bullet   
\end{array} $ }   \]
is involutive with characters $(0,0,2)$ and we may consider again the system:   \\
\[  \left \{ \begin{array}{lcl}
  y_{33} - x^2y_1 &  =  & v  \\
  y_{12}  &  =  &  u  \\
  y_{11}  &  =  & w= u_{33} -v_{12}- x^2u_1 
  \end{array} \right.   \]
Instead of doing the same change of variables, writing out the system $R'_3$ and study its formal inegrability with corresponding $9 + 11=20$ CC for $(u,v,w)$, an elementary but tedious computation, we shall use a trick, knowing in advance that the generating CC {\it must} be of order $1+1=2$ because $g'_2$ had to get one prolongation in order to become involutive and thus $2$-acyclic.  \\  
  
 \noindent
$\bullet$ {\it Step} 2: It thus remains to find out the CC for $(u,v)$ in the initial inhomogeneous system. As we have used two prolongations in order to exhibit $R'_2$, we have second order formal derivatives of $u$ and $v$ in the right members. Now, from the above argument, we have second order CC for the new right members and could hope therefore for a fourth order generating CC. The trick is to use the three different brackets of operators that can be obtained. We have in a formal way:  \\
\[ \begin{array}{lcl}
[d_{33} - x^2 d_1,d_{12}] y  &  = &  y_{1233} - x^2y_{112} - y_{1233} + x^2y_{112}+y_{11}   \\
                                                 &  =  &  y_{11}  \\
                                                 &  =  &  u_{33} - x^2u_1 -v_{12}
\end{array}   \]

\[  [d_2,d_{11} ] y= 0 = d_2(d_{11}y )- d_{1}(d_{12}y)= (u_{233}-v_{122}-x^2u_{12} -u_1)-u_1  \] 
brings the {\it third order} CC:   \\
\[  A\equiv u_{233}-v_{122}- x^2u_{12} - 2u_1 =0  \]

\[  \begin{array}{lcl}
 [ d_{33} - x^2d_1, d_{11} ] y  &  =  & y_{1133} - x^2y_{111} -y_{1133}+x^2y_{111}  \\
                                                         &  =  &  0  \\
                                                        &  =  &  (d_{33} - x^2d_1)w - v_{11}  \\
                                                        &  =  &  u_{3333} - v_{1233} 2 x^2u_{133}+ x^2v_{112}- v_{11} +
       (x^2)^2u_{11}
\end{array}     \]
brings the {\it fourth order} CC:   \\
\[      B  \equiv  u_{3333} - v_{1233} - 2 x^2u_{133}+ x^2v_{112}- v_{11} +(x^2)^2u_{11}=0  \]
We have indeed the identity $A_{33} - x^2 A_1 - B_2=0$ and thus $(A,B)$ are differentially dependent, that is $B$ is a new generating fourth order CC which is not a consequence of the prolongations of $A$. Again, the total number of generating CC, that is $1+1=2$, has nothing to do with the differential transcendence degree of the CC differential module which is $dim(F_0)-dim(E)=2-1=1$.  \\  \\

\noindent
{\bf EXAMLE 2.3}: With the same $ n=3,m=1, q=2$ and $K=\mathbb{Q}(x^1,x^2,x^3)=\mathbb{Q}(x)$, we now prove that a slight change of the equations may provide quite important changes in the number and order of the CC. Such an example is the only one that we could have found in more than $40$ years of computing CC in mathematics and applications. For this, let us consider the new system:  \\
\[        y_{33} - x^2y_1=v,  \,\,\,\,\,  y_{22}=u  \]
Before starting, we first notice that it is a prioiri not evident to discover that $R=R_{\infty}$ is a finite dimensional vector space over $K$ with $dim_K(R)=6$. However such a result can be obtained by direct integration (Compare to the Janet example treated in the introduction of [K1]).  \\

\noindent
$\bullet$ {\it Step} 1: The symbol $g_2$ is defined by $v_{33}=0, v_{22}=0$ may not be involutive {\it or} the coordinate system may not be $\delta$-regular. However, we obtain the Janet tabular for $g_2$:  \\
\[ \left \{ \begin{array}{lcl}
  v_{33} &  =  & 0  \\
  v_{22}  &  =  &  0 \\
  \end{array} \right. \fbox{ $\begin{array}{lll}
  1 & 2 & 3  \\
  1 & 2 & \bullet  
  \end{array} $}   \]
  and thus the Janet tabular for $g_3$:  \\
  \[ \left \{ \begin{array}{lcl}
  v_{333} &  =  & 0  \\
  v_{233}  & = & 0   \\
  v_{223}  &  =  &  0 \\
  v_{222}   &  =  &  0  \\
  v_{133}  &  =  &  0  \\
  v_{122}   &  =  &  0
  \end{array} \right. \fbox{ $\begin{array}{ccc}
  1 & 2 & 3  \\
  1  &  2  &  \bullet  \\
  1  &  2  &  \bullet  \\
  1 & 2 & \bullet   \\
  1  &  \bullet  &  \bullet  \\
  1  &  \bullet  &  \bullet  \\
  \end{array} $ }   \] 
 We let the reader check as an exercise that $g_2$ is not $2$-acyclic by counting the dimensions in the long sequence: \\
 \[  0 \rightarrow g_4  \stackrel{\delta}{\rightarrow} T^* \otimes g_3 \stackrel{\delta}{\rightarrow} {\wedge}^2T^*\otimes g_2 \stackrel{\delta}{\rightarrow} {\wedge}^3T^*\otimes T^*\]
and that $g_3$ is involutive, thus $2$-acyclic, with characters $(0,0,4)$ as in the previous example. It follows that $dim(g_2)=6-2=4$, $dim(g_3)=0+0+4=4, dim(g_4)=4, ...$. We obtain from the main theorem ${\rho}_r(R^{(1)}_2)= R^{(1)}_{2+r}$. It is easy to chek that $R^{(1)}_2=R_2$ with $dim(R_2)=8$, $dim(R_3)=8+4=12$, $dim(R_4)=11+4=15$, $dim(R_5)=56-39=13+4=17$, ... . We have after two prolongation:   \\
 \[ \left \{ \begin{array}{lcl}
  \,\,\,\,\,y_{2233}  & = &  u_{33}     \\
- y_{2233} + x^2y_{122} + 2y_{12} & = & - v_{22}  \\
-x^2y_{122}  &  =  &  - x^2u_1 
  \end{array}\right. \]
 and thus $2y_{12}=u_{33}- v_{22} - x^2u_1$. We may thus consider the new second order system $
 R'_2=R^{(2)}_2\subset R_2$ with a strict inclusion and $dim(R'_2)=7$:    \\
\[        y_{33} - x^2y_1=v,  \,\,\,\,\,  y_{22}=u , \,\,\,\, y_{12}=w, \,\,\,\,  2w=u_{33}-v_{22} - x^2 u_1 \] 
 
We may start again with $R'_2$ and study its symbol $g'_2 $ defined by the 3 linear equations with the following Janet tabular obtained after doing the same change of local coordinates as before:  \\
\[  \left \{ \begin{array}{lcl}
  v_{33} &  =  & 0  \\
  v_{22}  &  =  &  0 \\
  v_{12} &  =  &  0
  \end{array} \right. \fbox{ $\begin{array}{ccc}
  1 & 2 & 3  \\
  1 & 2 & \bullet  \\
  1  &  \bullet & \bullet   
\end{array} $ }   \]
This symbol is not involutive but its prolongation $g'_3$, defined by the $8$ equatios:  \\
\[  \left \{ \begin{array}{lcl}
  v_{333} &  =  & 0  \\
  v_{233}  & =  & 0   \\
  v_{223}  &  =  &  0 \\
  v_{222} &  =  &  0  \\
  v_{133} & = & 0  \\
  v_{123}   & = & 0  \\
  v_{122}  & = & 0  \\
  v_{112}  & = & 0
  \end{array} \right. \fbox{ $\begin{array}{ccc}
  1 & 2 & 3  \\
  1 & 2 & \bullet  \\
  1 & 2 & \bullet \\
  1 & 2 & \bullet  \\
  1 & \bullet & \bullet  \\
  1 & \bullet & \bullet  \\
  1 & \bullet & \bullet  \\
  1  &  \bullet & \bullet   
\end{array} $ }   \]
is involutive with characters $(0,0,2)$ and we may consider again the system:   \\
\[  \left \{ \begin{array}{lcl}
  y_{33} - x^2y_1 &  =  & v  \\
  y_{22}  &  =  &  u  \\
  y_{12}  &  =  & w
  \end{array} \right.   \]
with $2w=u_{33}-v_{22} -x^2u_1$. As before, instead of writing out the system $R'_3$ and studying its formal inegrability by an elementary but tedious computation, we shall use a trick, knowing in advance that the generating CC {\it must} be of order at least $1+1=2$ because $g'_2$ had to get one prolongation in order to become involutive and thus $2$-acyclic.  \\

 \noindent
$\bullet$ {\it Step} 2: It thus remains to find out the CC for $(u,v)$ in the initial inhomogeneous system. As we have used two prolongations in order to exhibit $R'_2$, we have second order formal derivatives of $u$ and $v$ in the right members. Now, from the above argument, we have second order CC for the new right members and could hope therefore for a fourth order generating CC. The trick is to use the three different brackets of operators that can be obtained. We obtain in a formal way:  \\
\[ \begin{array}{lcl}
[d_{33} - x^2 d_1,d_{22}] y  &  = &  y_{2233} - x^2y_{122} - y_{2233} + x^2y_{122}+2y_{12}   \\
                                                 &  =  &  2y_{12}  \\
                                                 &  =  &  u_{33} - x^2u_1 -v_{22}
\end{array}   \]
Then:   \\
\[ 2 [d_2,d_{12} ] y= 0 = d_2(2d_{12}y )- 2d_{1}(d_{22}y)= (u_{233}-v_{222}-x^2u_{12} -u_1)-2u_1  \] 
brings the {\it third order} CC:   \\
\[  A\equiv u_{233}-v_{222}- x^2u_{12} - 3u_1 =0  \]

\[  \begin{array}{lcl}
 2 [ d_{33} - x^2d_1, d_{12} ] y  &  =  & 2 y_{1233} - 2 x^2y_{112} - 2y_{1233}+ 2 x^2y_{112}  + 2y_{11} \\ 
                                                     &  =  &  2y_{11}  \\
                                                        &  =  &  2(d_{33} - x^2d_1)w - 2v_{12}  \\
                                                        &  =  &  u_{3333} - v_{2233} - 2 x^2u_{133}+ x^2v_{122}- 2v_{12} +
       (x^2)^2u_{11}
\end{array}     \]
brings the new first order equation:  
\[      2 y_{11}= 2 \omega = u_{3333} - v_{2233} - 2 x^2u_{133}+ x^2v_{122}- 2v_{12} +(x^2)^2u_{11}=0  \]
Accordingly, we may start afresh with the new system $R"_2= R^{(4)}_2\subset R'_2 \subset R_2$ which is {\it surprisingly} of finite type with $dim(R"_2)= 10-4=6$, $g"_3=0$ and defined by the $4$ second order PD equations:   \\
\[     y_{33} - x^2 y_1 =v, \,\,\,\,  y_{22} = u, \,\,\,\, y_{12}=w, \,\,\,\,  y_{11}=\omega , \,\,\, w\in j_2(u,v), \,\,\, \omega \in j_4(u,v)    \]
We obtain therefore: \\
\[   \begin{array}{rcl}
2 y_{112}  &  =  &  u_{23333} - v_{22233} - 2x^2u_{1233} - 2u_{133}+ x^2v_{1222} - v_{122} + (x^2)^2u_{112} + 2x^2 u_{11}     \\
- 2 y_{112} & =  & - u_{133} +v_{122} + x^2 u_{11}
\end{array}   \]
and thus a CC of order $5$, namely:  \\
\[ B \equiv u_{23333} - v_{22233} - 2 x^2u_{1233}- 3u_{133} + x^2v_{1222} +(x^2)^2u_{112}+
3x^2u_{11} = 0   \]
However, we have indeed the identity $A_{33} - x^2 A_1 - B=0$ and thus $B\in j_2(A)$, that is $B$ is a {\it not} a new generating fifth order CC as it is only a consequence of the prolongations of $A$.
Using now the bracket of operators $[d_{11}, d_{33}] = 0$ that has not been already used, we get:   \\
\[  \begin{array}{lcl}
2y_{1133} & = & u_{333333} - v_{223333} - 2x^2u_{13333} + x^2v_{12233} - 2v_{1233} + (x^2)^2u_{1133} \\
- 2 y_{1133} + 2 x^2 y_{111} & = & - 2 v_{11}   \\
-2 x^2y_{111} & = &  -x^2 u_{13333} + x^2 v_{12233} + 2(x^2)^2u_{1133} - (x^2)^2v_{1122} + 2 x^2v_{112} - (x^2)^3u_{111}
\end{array}    \]
We obtain therefore a new sixth order CC:   \\
\[   \begin{array}{cl}
C  \equiv & u_{333333} -v_{223333}   \\
              &  - 3x^2u_{13333} + 2 x^2v_{12233}    \\
              &  - 2v_{1233} + 3 (x^2)^2u_{1133} - (x^2)^2v_{1122}   \\
              &  + 2 x^2v_{112} - (x^2)^3u_{111}  \\
              &  - 2v_{11}    \\
              &   =0
              \end{array}   \]
which cannot be a differential consequence of $A$. After tedious computations, one can find the differential identity:  \\
\[     A_{3333} - 2 x^2 A_{133} + (x^2)^2 A_{11} - C_2 =0  \]

The corresponding simplest free resolution, written with differential modules, is thus:   \\
\[  0 \rightarrow D \underset 4 {\longrightarrow} D^2  \underset 6{\longrightarrow} D^2 \underset 2{\longrightarrow} D \stackrel{p}{\longrightarrow} M \rightarrow  0   \]

Again, the total number of generating CC, that is $1+1=2$, has nothing to do with the differential transcendence degree of the CC differential module which is still $dim(F_0)-dim(E)=2-1=1$ because $rk_D(M)=0$.  \\  \\  \\

\newpage

\noindent
{\bf 3) MATHEMATICAL TOOLS} \\

Instead of starting with a linear system $R_q\subset J_q(E)$ of order $q$ on $E$, let us start with a bundle map $\Phi: J_q(E)  \rightarrow F_0$ with 
$R_q=ker(\Phi)$ , let us consider the linear PD operator  $  {\cal{D}}: E \stackrel{j_q}{\longrightarrow}J_q(E) \stackrel{\Phi}{\longrightarrow} F_0$. Using the canonical inclusion $J_{q+r)(E} \subset J_r(J_q(E) )$, let us define the $r$-{\it prolongation} ${\rho}_r(\Phi):J_{q+r}(E) \rightarrow J_r(J_q(E)) \stackrel{J_r(\Phi)}{\longrightarrow} J_r(F_0)$. The general case of the successive prolongations with $r\geq0$ is described by the following commutative and exact diagram:  \\
\[  \begin{array}{rcccccccl}
    &  0  &  &  0  &  &  0  &  &  &   \\
     &  \downarrow &  & \downarrow &  &  \downarrow  &  &   &     \\
0 \rightarrow  &  g_{q+r+1}  & \rightarrow & S_{q+r+1}T^*\otimes E & \stackrel{{\sigma}_{r+1}(\Phi)}{\longrightarrow} &  S_{r+1}T^*\otimes F_0 & \rightarrow &   h_{r+1} &  \rightarrow 0 \\
&  \downarrow &  & \downarrow &  &  \downarrow  &  &  \downarrow    &     \\
0 \rightarrow &  R_{q+r+1} & \rightarrow  & J_{q+r+1}(E) & \stackrel{{\rho}_{r+1}(\Phi)}{\longrightarrow } &  J_{r+1}(F_0) &  \rightarrow &  Q_{r+1}  &  \rightarrow 0  \\
    &  \downarrow &  & \hspace{12mm}\downarrow {\pi}^{q+r+1}_{q+r}&  & \hspace{8mm} \downarrow  {\pi}^{r+1}_{r}  &  & \downarrow   &     \\
0 \rightarrow &  R_{q+r} & \rightarrow  & J_{q+r}(E) & \stackrel{{\rho}_r(\Phi)}{\longrightarrow}  &  J_r(F_0) &  \rightarrow &  Q_{r}  & \rightarrow 0 \\ 
   &  &  & \downarrow &  &  \downarrow  &  & \downarrow  & \\
   &  &  &  0  && 0  && 0  &
\end{array}   \]
with {\it symbol}-map induced in the upper {\it symbol}-sequence.  \\
Chasing in this diagram while applying the "{\it snake}" lemma ([11],[14],[22]), we obtain the {\it long exact connecting sequence}:  \\
 \[  0  \rightarrow  g_{q+r+1}  \rightarrow R_{q+r+1}  \rightarrow R_{q+r}  \rightarrow h_{r+1} \rightarrow Q_{r+1} \rightarrow Q_r \rightarrow 0  \]
which is thus connecting in a tricky way FI ({\it lower left}) with CC ({\it upper right}). Needless to say that absolutely no classical procedure can produce such a result which is thus totally absent from the GR papers already quoted. \\  
Setting $H(R_{q+r})=R_{q+r}/{\pi}^{q+r+1}_{q+r}(R_{q+r+1})$, we have equivalently the shorter long exact sequence:  \\
\[  0 \rightarrow H(R_{q+r}) \rightarrow h_{r+1} \rightarrow Q_{r+1}  \rightarrow Q_r \rightarrow 0  \]
As a posible interpretation, $dim(Q_r)$ is the total number of CC of order $0, 1, ...$ up to $r$ included. However, the problem to solve is to study the structure of the projective limit of vector bundles made by the induced epimorphisms $Q_{r+1} \rightarrow Q_r$. Of course, as it is mostly realized in the examples, we have to suppose that $R_q$ is {\it sufficiently regular} in such a way that the $R_{q+r}$ are vector bundles $\forall r\geq 0${\it and} that the  $R^{(s)}_{q+r}= {\pi}^{q+r+s}_{q+r}(R_{q+r+s}$ are also vector bundles, such a situation being in particular always realized when $R_q\subset J_q(E)$ or ${\cal{D}}$ are defined over a differential field $K$. In this case, introducing the filtred noetherian ring $D=K[d_1,...,d_n]=K[d]$ of differential operators with coefficients in $K$, we may introduce a differential module $M$ with induced filtration  $0 =M_0\subseteq M_1 \subseteq ... \subseteq M_q \subseteq ... \subseteq m_{\infty}=M$ in such a way that the system $R=R_{\infty}=hom_K(M,K)$ associated with $M$ with $R_q =hom_K(M_q,K)$ is of course {\it automatically} FI ({\it care}). Following Macaulay in ([10]), we have already proved in many places ([14],[19]) that $R$ is a differential module for the {\it Spencer operator} $d:R \rightarrow T^*\otimes R:f \rightarrow dx^i\otimes d_if$ with $d_i:R \rightarrow R:R_{q+1} \rightarrow R_q$ defined by:  \\
\[   (d_if)^k_{\mu}={\partial}_if^k_{\mu}- f^k_{\mu +1_i} \in K   \]
It is important to notice that such an operator/system is far from being formally integrable because:   \\
\[  {\partial}_i({\partial}_jf^k_{\mu} -f^k_{\mu +1_j}) - {\partial}_j ( {\partial}_if^k_{\mu} - f^k_{\mu +1_i})=
{\partial}_jf^k_{\mu +1_i} -{\partial}_if^k_{\mu +1_j}   \]
As can be seen from the examples previously presented, starting with ${\Psi}_r$ for a given $r$, the main problem is to compare the epimorphism 
${\Psi}_{r+1}:J_{r+1}(F_0) \rightarrow Q_{r+1}$ with the morphism ${\rho}_1({\Psi}_r): J_{r+1}(F_0) \rightarrow J_1(Q_r)$ in the following commutative diagram {\it which may not be exact}:  \\
\[  \begin{array}{rcccccccl}
    &  0  &  &  0  &  &  0  &  & 0  &   \\
     &  \downarrow &  & \downarrow &  &  \downarrow  &  & \downarrow   &     \\
0 \rightarrow  &  g_{q+r+1}  & \rightarrow & S_{q+r+1}T^*\otimes E & \stackrel{{\sigma}_{r+1}(\Phi)}{\longrightarrow} &  S_{r+1}T^*\otimes F_0 & \stackrel{{\sigma}_1({\Psi}_r)}{\longrightarrow} &   T^*\otimes Q_r &   \\
&  \downarrow &  & \downarrow &  &  \downarrow  &  &  \downarrow    &     \\
0 \rightarrow &  R_{q+r+1} & \rightarrow  & J_{q+r+1}(E) & \stackrel{{\rho}_{r+1}(\Phi)}{\longrightarrow } &  J_{r+1}(F_0) & 
 \stackrel{{\rho}_1({\Psi}_r)}{\longrightarrow} &  
J_1(Q_r)  &    \\
    &  \downarrow &  & \hspace{12mm}\downarrow {\pi}^{q+r+1}_{q+r}&  & \hspace{8mm} \downarrow  {\pi}^{r+1}_r&  & \downarrow   &     \\
0 \rightarrow &  R_{q+r} & \rightarrow  & J_{q+r}(E) & \stackrel{{\rho}_r(\Phi)}{\longrightarrow}  &  J_r(F_0) &  \stackrel{{\Psi}_r}{\longrightarrow} &  Q_{r}  & \rightarrow 0 \\ 
   &  &  & \downarrow &  &  \downarrow  &  & \downarrow  & \\
   &  &  &  0  && 0  && 0  &
\end{array}   \]
where the central row is induced from the long exact sequence:  \\
\[  0 \rightarrow J_1(R_{q+r}) \rightarrow  J_1(J_{q+r}(E)) \rightarrow J_1(J_r(F_0))  \rightarrow J_1(Q_r) \rightarrow  0  \]
and may not be exact.  \\

 \noindent
 {\bf PROPOSITION 3.1}: We have only in general:   \\
\[  B_{r+1} = im({\rho}_{r+1}(\Phi))  = ker({\Psi}_{r+1})\subseteq ker({\rho}_1({\Psi}_r) ) = {\rho}_1(B_r) \]

\noindent
{\it Proof}: Denoting the Spencer operator by $d$ in place of the standard notation $D$ of the literature that could be confused with the ring $D$ of differential operators, we have the following commutative diagram:  \\
\[  \begin{array}{rccccccl}
   
0 \rightarrow &  R_{q+r+1} & \rightarrow  & J_{q+r+1}(E) & \stackrel{{\rho}_{r+1}(\Phi)}{\longrightarrow } & B_{r+1}& \subset &  J_{r+1}(F_0) \\
    & \hspace{2mm} \downarrow  d  &  & \hspace{2mm}\downarrow  d  &  & \hspace{2mm} \downarrow  d  &  & \hspace{2mm}\downarrow d  \\
0 \rightarrow &  T^* \otimes R_{q+r} & \rightarrow  & T^* \otimes J_{q+r}(E) & \stackrel{{\rho}_r(\Phi)}{\longrightarrow}  & T^*\otimes B_r & \subset & T^* \otimes J_r(F_0)     
\end{array}   \]
As $B_{r+1}$ projects {\it onto} $B_r$ and $d_iB_{r+1} \subset B_r$, it follows from ([13], Propositions $10$, p $83$) or ([14], Remark 2.9, p 315)  that Ê$B_{r+1}\subseteq {\rho}_1(B_r)$. We have thus a projective limit of systems, each one being defined by more equations than the preceding oneÊÊand such a procedure must finish with a FI system that can even be prolongued, as we shall see in the examples, in order to obtain an involutive systemthat may be used to start a Janet sequence. The decision to stop is provided by the maximum order of the CC obtained, namely of order bounded by $r+s+1=t$ if the system $R^{(s)}_{q+r}$ is involutive or at least with a $2$-acyclic symbol. \\ÊÊÊÊÊÊÊ
\hspace*{12cm}  \fbox{}   \\
ÊÊÊÊÊÊÊÊÊÊÊÊÊÊÊÊÊÊÊÊÊÊÊÊÊÊÊÊÊÊÊÊÊÊÊÊÊÊÊÊÊÊÊÊÊÊÊÊÊÊÊÊÊÊÊÊÊÊÊÊÊÊÊÊÊÊÊÊÊÊÊÊÊÊÊÊÊÊÊÊÊÊÊÊÊÊÊÊÊÊÊÊÊÊÊÊÊÊÊÊÊÊÊÊÊÊÊÊÊÊÊÊÊÊÊÊÊÊÊÊÊÊÊÊÊÊÊÊÊÊÊÊÊÊÊÊÊÊÊÊÊÊÊÊÊÊÊÊÊÊÊÊÊÊÊÊÊÊÊÊÊÊÊÊÊÊÊÊÊÊÊÊÊÊÊÊÊÊÊÊÊÊÊÊÊÊÊÊÊÊÊÊÊÊÊÊÊÊÊÊÊÊÊÊÊÊÊÊÊÊÊÊÊÊÊÊÊÊÊ

The idea is to use the {\it composite morphism} $S_rT^*\otimes F_0 \rightarrow Q_r$  while chasing in order to prove that any element of $J_{r+1}(F_0)$ killed by ${\rho}_1({\Psi}_r)$ can be decomposed into the sum of an element in $im({\rho}_r(\Phi))$ plus an element in $S_{r+1}T^*\otimes F_0$ killed by ${\sigma}_1({\Psi}_r)$. With more details, setting for simplicity $S_rT^*\otimes F_0=S_r(F_0)$, introducing the {\it coboundary} bundle $B(S_{r+1}(F_0))=im({\sigma}_{r+1}(\Phi))$ and the {\it cocycle} bundle $Z(S_{r+1}(F_0))=ker({\sigma}_1({\Psi}_r))$, we may define the corresonding {\it cohomology} bundle $H(S_{r+1}(F_0))=Z(S_{r+1}(F_0))/B(S_{r+1}(F_0))$. We may also define similarly $H(J_{r+1}(F_0))= ker({\rho}_1({\Psi}_r))/im({\rho}_{r+1}(\Phi))$ and we obtain the following crucial proposition (See [MSK], Example 2.A.9):  \\

\noindent
{\bf PROPOSITION 3.2}:  There exists a short exact sequence:   \\
\[ 0 \rightarrow H(R_{q+r}) \rightarrow  H(S_{r+1}(F_0)) \rightarrow H(J_{r+1}(F_0))  \rightarrow 0 \]

Le us now deal with the symbol cohomology by chasing in the following commutative diagram:  \\
 
\[  \begin{array}{rccccccc}
   &  0 & & 0 & & 0 &  &  0   \\
   & \downarrow & & \downarrow & & \downarrow & & \downarrow   \\
0\rightarrow & g_{q+r+1} & \rightarrow &  S_{q+r+1}T^*\otimes E & \rightarrow & S_{r+1}T^*\otimes F_0& \rightarrow & T^*\otimes Q_r    \\
   & \hspace{2mm}\downarrow  \delta  & & \hspace{2mm}\downarrow \delta & &\hspace{2mm} \downarrow \delta & \searrow &\parallel   \\
0\rightarrow& T^*\otimes g_{q+r}&\rightarrow &T^*\otimes S_{q+r}T^*\otimes E & \rightarrow &T^*\otimes S_rT^*\otimes F_0 &\rightarrow & T^*\otimes Q_r   \\
   &\hspace{2mm} \downarrow \delta &  &\hspace{2mm} \downarrow \delta & &\hspace{2mm}\downarrow \delta &  & \downarrow    \\
0\rightarrow & {\wedge}^2T^*\otimes g_{q+r-1} & \rightarrow & {\wedge}^2T^*\otimes S_{q+r-1}T^*\otimes E & \rightarrow & {\wedge}^2T^*\otimes S_{r-1}T^*\otimes  F_0 & & 0   \\
   &\hspace{2mm}\downarrow \delta  &  & \hspace{2mm} \downarrow \delta  &  &  & &    \\
0\rightarrow & {\wedge}^3T^*\otimes g_{q+r-2} & =  & {\wedge}^3T^*\otimes S_{q+r-2}T^*\otimes E  &   &   &  &      
\end{array}  \]
where neither the first nor the second upper columns may be exact and where the left column may not be
exact, unless $gq$ is involutive or $2$-acyclic.  Chasing with the same notations, we obtain:  \\

\noindent
{\bf PROPOSITION 3.3}: There exists an exact sequence:   \\
\[    0  \rightarrow H^2(g_{q+r-1}) \rightarrow H(S_{r+1}(F_0)) \rightarrow T^*\otimes H(S_r(F_0))     \]
The upper left arrows are not in general epimorphisms and it may be sometimes useful to consider $h_r$ as a kind of symbol in the more abstract diagram:   \\
 
{\small \[  \begin{array}{rcccccccl}
   &  0 & & 0 & & 0 &  &  &  \\
   & \downarrow & & \downarrow & & \downarrow & & &  \\
0\rightarrow & g_{q+r+1} & \rightarrow &  S_{q+r+1}T^*\otimes E & \rightarrow & S_{r+1}T^*\otimes F_0& \rightarrow & h_{r+1} & \rightarrow 0 \\
   & \hspace{2mm}\downarrow  \delta  & & \hspace{2mm}\downarrow \delta & &\hspace{2mm} \downarrow \delta &  &\hspace{2mm} \downarrow \delta &  \\
0\rightarrow& T^*\otimes g_{q+r}&\rightarrow &T^*\otimes S_{q+r}T^*\otimes E & \rightarrow &T^*\otimes S_rT^*\otimes F_0 &\rightarrow & T^*\otimes h_r  &\rightarrow 0 \\
   &\hspace{2mm} \downarrow \delta &  &\hspace{2mm} \downarrow \delta & &\hspace{2mm}\downarrow \delta &  &  \hspace{2mm} \downarrow \delta &  \\
0\rightarrow & {\wedge}^2T^*\otimes g_{q+r-1} & \rightarrow & {\wedge}^2T^*\otimes S_{q+r-1}T^*\otimes E & \rightarrow & {\wedge}^2T^*\otimes S_{r-1}T^*\otimes  F_0 & \rightarrow & {\wedge}^2T^*\otimes h_{r-1} &\rightarrow 0  \\
   &\hspace{2mm}\downarrow \delta  &  & \hspace{2mm} \downarrow \delta  &  &\hspace{2mm}\downarrow \delta   & & &  \\
0\rightarrow & {\wedge}^3T^*\otimes g_{q+r-2} & \rightarrow  & {\wedge}^3T^*\otimes S_{q+r-2}T^*\otimes E  & \rightarrow & {\wedge}^3T^*\otimes S_{r-2}T^*\otimes F_0  &   &  &      \\
   & \hspace{2mm}  \downarrow \delta & &  \hspace{2mm} \downarrow \delta &  &  &  &  & \\
0 \rightarrow &  {\wedge}^4T^*\otimes g_{q+r-3} & \rightarrow & {\wedge}^4T^*\otimes S_{q+r-3}T^*\otimes E  &   &  &  &  & 
\end{array}  \]   } 
where the rows are now exact. However, understanding the meaning of $h_r$ as a kind of new symbol may not be possible unless $h_{r+1} \rightarrow T^* \otimes h_r$ is a monomorphism, that is when $g_q$ is $2$-acyclic {\it and} $h_r$ is $1$-acyclic, that is when $g_q$ is {\it also} $3$-acyclic (or involutive). Once more, we understand the crucial importance of $2$-acyclicity but we recall that {\it the only symbol known to be $2$-acyclic without being involutive is the symbol of the conformal Killing system whenever $n \geq 4$ , which is also $3$-acyclic whenever} $n\geq 5$ ([12],[13],[17],[19]).  \\  \\

\noindent
{\bf 4) APPLICATIONS}  \\

\noindent
{\bf MACAULAY EXAMPLE REVISITED}:   \\

With $m=1, n=3, q=2, K= \mathbb{Q}$, let us let us introduce two operators $P,Q \in D=K[d_1,d_2,d_3]$ and consider the second order system $R_2\subset J_2(E)$ used by Macaulay as early as in $1916$ ([10]):   \\
\[  Qy \equiv y_{33} = v, \hspace{2cm}  Py\equiv y_{13} - y_2 =u \]
We have the strict inclusions:  \\
\[    R^{(2)}_2 \subset R^{(1)}_2 \subset R_2 \subset J_2(E)                                  \hspace{2cm}  6 < 7 < 8 < 10    \]
As $g^{(2)}_2 \subset  g^{(1)}_2 \subset g_2$ are involutive, we obtain ${\rho}_r(R^{(2)}_2)=R^{(2)}_{r+2}$ by using the {\it Prolongation/Projection} (PP) procedure. We exhibit the parametric jets of the bundles that will be used in the following diagrams:   \ \
\[  {par}_2=\{  y, y_1, y_2, y_3, y_{11}, y_{12}, y_{22}, y_{23}  \}   \]
\[  {par}_3= \{  y, y_1, y_2, y_3, y_{11}, y_{12}, y_{22}, y_{111}, y_{112}, y_{122}, y_{222}, y_{223}\}   \]
\[ {par}_4= \{  y, y_1, y_2, y_3, y_{11}, y_{12}, y_{111}, y_{112}, y_{122}, y_{222},y_{1111}, y_{1112}, y_{1122}, y_{1222}, y_{2222}, y_{2223}\}  \]
and thus $dim(R_2)=8, \,\, dim(R_3)=12, \,\, dim(R_4)= 16$. More generally, we let the reader prove that $dim(R_{r+2})=4r + 8, \forall r\geq 0$, thus $dim(R_{r+4})=4r+ 16,  \forall r\geq 0 $ and $dim(g_{r+4})=r+6, \forall r\geq 0$.  \\
We have the Janet tabular for $R^{(2)}_2$:  \\
\[  \left \{ \begin{array}{lcl}
 y_{33} & = & v \\
 y_{23}  &  = &  v_1 - u_3 \\
y_{22}& = &  v_{11}- u_{13} - u_2  \\
y_{13} - y_2  & = & u
 \end{array} \right.  \fbox { $ \begin{array}{ccc}
  1 & 2 & 3  \\
 1 & 2 & \bullet  \\
 1 & 2 & \bullet  \\
 1 & \bullet & \bullet 
\end{array} $ }   \]
The {\it two} CC are:  \\
\[  A\equiv v_{13} - v_2 - u_{33}=0, \hspace{2cm}  A_1 \equiv   v_{113} - v_{12} - u_{133}=0   \]
while the other ones are what we called {\it identity to zero} like:  \\
\[  d_2(y_{13} - y_2) -d_1(y_{23})+d_2y_2=0=u_2 -(v_{11} - u_{13}) +(v_{11} - u_{13} - u_2)=0   \]
There is thus only {\it one} generating CC of order $2$, namely $A=0$, given by the commutation relation $P \circ Q - Q \circ P \equiv 0 $ and the corresponding operator ${\cal{D}}_1$ is thus surely formally surjective. Setting $F_1=Q_2$, we obtain the following diagram with exact central and lower rows whenever $r\geq 1$.      \\
{\small  \[  \begin{array}{rcccccccl}
    &  0  &  &  0  &  &  0  &  & 0  &   \\
     &  \downarrow &  & \downarrow &  &  \downarrow  &  &\downarrow   &     \\
0 \rightarrow  &  g_{r+4}  & \rightarrow & S_{r+4}T^*\otimes E & \rightarrow &  \fbox{$S_{r+2}T^*\otimes F_0$} & \rightarrow &  S_rT^*\otimes F_1 &  \rightarrow 0 \\
&  \downarrow &  & \downarrow &  &  \downarrow  &  &  \downarrow    &     \\
0 \rightarrow &  R_{r+4} & \rightarrow  & J_{r+4}(E) & \rightarrow  &  J_{r+2}(F_0) &  \rightarrow &  J_r(F_1)  &  \rightarrow 0  \\
    &  \downarrow &  & \downarrow &  &  \downarrow  &  & \downarrow   &     \\
0 \rightarrow &  R_{r+3} & \rightarrow  & J_{r+3}(E) & \rightarrow  &  J_{r+1}(F_0) &  \rightarrow &  J_{r-1}(F_1)  &  \rightarrow 0  \\ 
   &  &  & \downarrow &  &  \downarrow  &  &  \downarrow   & \\
   &  &  &  0  && 0  && 0 &
\end{array}   \]}
{\scriptsize  \[  \begin{array}{rcccccccl}
    &  0  &  &  0  &  &  0  &  &  0 &   \\
     &  \downarrow &  & \downarrow &  &  \downarrow  &  & \downarrow  &     \\
0 \rightarrow  &  r+6  & \rightarrow &(r+5)(r+6)/2 & \rightarrow & \fbox{(r+3)(r+4)} & \rightarrow &  (r+1)(r+2)/2&  \rightarrow 0 \\
&  \downarrow &  & \downarrow &  &  \downarrow  &  &  \downarrow    &     \\
0 \rightarrow &  4r+16 & \rightarrow  & (r+5)(r+6)(r+7)/6  & \rightarrow  & (r+3)(r+4)(r+5)/3  &  \rightarrow & (r+1)(r+2)(r+3)/6  &  \rightarrow 0  \\
    &  \downarrow &  & \downarrow &  &  \downarrow  &  & \downarrow   &     \\
0 \rightarrow & 4r+12 & \rightarrow  &  (r+4)(r+5)(r+6)/6  & \rightarrow  & (r+2)(r+3)(r+4)/3  &  \rightarrow &  r(r+1)(r+2)/6  & \rightarrow  0   \\ 
   & \downarrow &  & \downarrow &  &  \downarrow  &  & \downarrow  & \\
   & \fbox{r+2} &  &  0  &  & 0  &  & 0  &  \\
    & \downarrow &  &    &  &   &  &  &  \\
     &  0  &  &    &  &   &  &  &  \\
\end{array}   \]}
{\scriptsize \[  \begin{array}{rcccccccl}
 & 0\hspace{2mm}& & 0 \hspace{2mm} && 0\hspace{2mm}& & 0 \hspace{2mm} & \\
 & \downarrow \hspace{2mm}& & \downarrow \hspace{2mm}& & \downarrow \hspace{2mm} & & \downarrow  \hspace{2mm}&\\
0 \rightarrow & g_{r+4} & \rightarrow & S_{r+4}T^*\otimes E & \rightarrow  &\fbox{$ S_{r+2}T^* \otimes F_0$} & \rightarrow  & S_rT^*\otimes F_1 & \rightarrow 0  \\
  &  \downarrow \delta  & & \downarrow \delta & & \downarrow  \delta & & \downarrow \delta  &\\
0 \rightarrow & T^*\otimes g_{r+3} & \rightarrow & T^* \otimes S_{r+3}T^* \otimes E & \rightarrow & T^*\otimes S_{r+1}T^*\otimes F_0 &\rightarrow & T^*\otimes S_{r-1}T^*\otimes F_1 & \rightarrow 0  \\
   &  \downarrow \delta  & & \downarrow \delta & & \downarrow  \delta &   &   \\
0 \rightarrow & \fbox{${\wedge}^2T^*\otimes g_{r+2}$} & \rightarrow &{\wedge}^2T^*\otimes S_{r+2}T^* \otimes E & \rightarrow & {\wedge}^2T^*\otimes S_rT^*\otimes F_0 &  &   & \\
  & \downarrow \delta &  & \downarrow \delta & & & &  & \\
 0 \rightarrow  & {\wedge}^3T^*\otimes  g_{r+1} & = &   {\wedge}^3T^*\otimes S_{r+1}T^*\otimes E & && & &  \\
  &\downarrow \hspace{2mm} & &\downarrow \hspace{2mm} & & & & &  \\
  & 0 \hspace{2mm} &  &  0 \hspace{2mm} &  &  &  &  &
\end{array}   \] }

\noindent
where $S_4T^*\otimes E \simeq S_4T^*$ and $F_1\simeq Q_2$ while $Q_1=0$ as there is no CC of order $1$. From the snake lemma and a chase, we obtain the {\it long exact connecting sequence} when $r=0$:  \\
\[  0  \rightarrow  g_4  \rightarrow R_4  \rightarrow R_3  \rightarrow h_2 \rightarrow F_1  \rightarrow 0  \]
\[  0  \rightarrow  6  \rightarrow 16 \rightarrow 12 \rightarrow 3 \rightarrow  1  \rightarrow 0  \]
relating FI ({\it lower left}) to CC ({\it upper right}). By composing the epimorphism $S_2T^* \otimes F_0 \rightarrow h_2$ with the epimorphism $h_2\rightarrow F_1$, we obtain an epimorphism $S_2T^*\otimes F_0 \rightarrow F_1$ and the long exact sequence:  \\
\[ 0 \rightarrow g_4 \rightarrow S_4T^*\otimes E \rightarrow  S_2T^*\otimes F_0 \rightarrow F_1 \rightarrow  0  \]  \\
which is nevertheless {\it not} a long ker/coker exact sequence by counting the dimensions as we have $6 - 15 + 12 -1= 2 \neq 0$.  \\
The above diagrams illustrate perfectly the three propositions of Section $2$. We have in particular:  \\
\[     H(J_{r+2}(F_0))=0 \,\,\, \Rightarrow \,\,\, H(J_{r+1}(F_0)=0, \hspace{1cm}  H(S_{r+2}(F_0))=H(R_{r+3})\neq 0   \]
and the formally exact sequence, which is nevertheless not strictly exact though $1-2+1=0$:  \\
\[  0  \rightarrow  \Theta \rightarrow E \underset 2{\stackrel{{\cal{D}}}{\longrightarrow}} F_0 \underset 2 {\stackrel{{\cal{D}}_1}{\longrightarrow}} 
F_1  \rightarrow 0  \]
\[  0  \rightarrow  \Theta \rightarrow 1\longrightarrow  2 \longrightarrow  1 \rightarrow 0  \]  \\
We remind the reader that, contrary to the situation met for FI systems where the exactness on the jet level is obtained inductively from the exactness on the symbol level, here we discover that {\it we may have the exactness on the jet level without having exactness on the symbol level}.  \\  \\

\noindent
{\bf EXAMPLE 2.1 REVISITED}: \\

First of all, let us compute the dimensions and the parametric jets that will be used in the following diagrams.
\[  par_1=par_2==\{ y, y_1, y_2, y_3, y^2_1\}, \,\,\,\,\,  \]
\[  n=dim(X)=3, \,\,\, m=dim (E)=2, \,\,\,dim(R_1)=dim(R_2)=3,\,\,\,  dim(g_1)=1, \,\,\, g_2=0 \Rightarrow g_3=0 \]

 \[  \begin{array}{rccccccccccl}
 &&&& 0 &&0&&0&  & \\
 &&&& \downarrow && \downarrow && \downarrow  & & \\
  && 0 & \longrightarrow &S_3T^*\otimes E&\longrightarrow&S_2T^*\otimes F_0 &\longrightarrow & T^*\otimes Q_1 & \longrightarrow 0      \\
   &&\downarrow && \downarrow & & \downarrow &\searrow & \downarrow & &  & \\
0& \longrightarrow &R_3&\longrightarrow &J_3(E)& \longrightarrow &J_2(F_0) &\longrightarrow & J_1(Q_1) & \longrightarrow 0  \\
    & & \downarrow && \hspace{5mm}\downarrow {\pi}^3_4 & &\hspace{5mm} \downarrow {\pi}^2_1 & & \downarrow  &  \\
   0 & \longrightarrow & R_2& \longrightarrow & J_2(E) & \longrightarrow &J_1(F_0)& \longrightarrow &Q_1& \longrightarrow  0 \\
   && \downarrow && \downarrow & & \downarrow & & \downarrow &      \\
   &&0&& 0 && 0 && 0  &  
   \end{array}     \]
 \[  \begin{array}{rccccccccccl}
 &&&& 0 &&0&& 0 &  & \\
 & &&& \downarrow && \downarrow &&\downarrow &  & \\
  && 0 & \longrightarrow &20 &\longrightarrow& 30 &\longrightarrow & 9 &  \longrightarrow 0      \\
   &&\downarrow && \downarrow & & \downarrow & \searrow & \downarrow & &   & \\
0& \longrightarrow &3&\longrightarrow & 40 & \longrightarrow & 50 &\longrightarrow & 12 &   \longrightarrow 0 \\
    & & \downarrow && \hspace{5mm}\downarrow {\pi}^3_2 & &\hspace{5mm} \downarrow {\pi}^2_1 & & \downarrow  &  \\
   0 & \longrightarrow & 3 & \longrightarrow & 20 & \longrightarrow & 20 & \longrightarrow &  3  & \longrightarrow  0 \\
   &&\downarrow && \downarrow & & \downarrow & & \downarrow &      \\
   && 0 && 0 && 0 && 0  &  
   \end{array}     \]

 \[  \begin{array}{rccccccccccl}
 &&&& 0 &&0&& 0 &  & \\
 &&&& \downarrow && \downarrow && \downarrow& & \\
  && 0 & \rightarrow &S_3T^*\otimes E&\rightarrow&S_2T^*\otimes F_0 &\rightarrow & T^*\otimes Q_1& \rightarrow 0
    \\
   &&&& \hspace{3mm} \downarrow \delta & & \hspace{3mm} \downarrow \delta & & \parallel & &  & \\
&  &0&\rightarrow &T^*\otimes S_2T^*\otimes E & \rightarrow &T^*\otimes T^* \otimes F_0 &\rightarrow & T^*\otimes Q_1 &   \rightarrow 0 \\
    & & && \hspace{3mm}\downarrow \delta  & &\hspace{3mm} \downarrow \delta& & \downarrow  &  \\
   0 & \rightarrow & {\wedge}^2T^*\otimes g_1& \rightarrow & {\wedge}^2T^*\otimes T^* \otimes E & \rightarrow & {\wedge}^2T^*\otimes F_0  & \rightarrow & 0 &  \\
   && \hspace{3mm} \downarrow \delta && \hspace{3mm} \downarrow \delta & & \downarrow & & &   &   \\
 0   & \rightarrow & {\wedge}^3T^*\otimes E  & = & {\wedge}^3T^*\otimes E & \rightarrow & 0  & &&&   \\
       &   &  \downarrow &  & \downarrow  &  &  &  &  &  &  \\
   &&0&& 0 && && &  &
   \end{array}     \]

 \[  \begin{array}{rccccccccccl}
 &&&& 0 &&0&&0&  & \\
 &&&& \downarrow && \downarrow &&\downarrow & & \\
  && 0 & \rightarrow &20&\rightarrow&\fbox{30} &\rightarrow & 9 &
   \rightarrow  0  \\
   &&&& \hspace{3mm} \downarrow \delta & & \hspace{3mm} \downarrow \delta & & \parallel& &    & \\
&  &0&\rightarrow &36 & \rightarrow &45 &\rightarrow & 9 &   \rightarrow 0 \\
    & & && \hspace{3mm}\downarrow \delta  & &\hspace{3mm} \downarrow \delta& & \downarrow  &  \\
   0 & \rightarrow & \fbox{3} & \rightarrow & 18 & \rightarrow & 15 & \rightarrow & 0 &  \\
   && \hspace{3mm} \downarrow \delta && \hspace{3mm} \downarrow \delta & & \downarrow & & &   &   \\
 0   & \rightarrow &  2  & = &  2  & \rightarrow & 0  & &&&   \\
       &   &  \downarrow &  & \downarrow  &  &  &  &  &  &  \\
   &&0&& 0 && && &  &
   \end{array}     \]
It is not at all evident to study these diagrams. We have $dim(B(S_2(F_0))=20 < dim(Z(S_2(F_0))= 30-9=21$.  We have already proved that $ dim (H(S_2(F_0))=21-20=1=dim(H^2(g_1))=3-2=1$, a result not evident at first sight explaining why the only second order additional generating CC is nothing else than the Riemann tensor in dimension equal to $2$. \\

We have explained in ([20]) that such a system has its origin in the study of the integration of the Killing system for the Schwarzschild metric, which is not FI. With more details, let us use the Boyer-Lindquist coordinates $(t,r,\theta, \phi)=(x^0,x^1,x^2,x^3)$ instead of the Cartesian coordinates $(t,x,y,z)$ and consider the Schwarzschild metric $\omega= A(r)dt^2-(1/A(r))dr^2 - r^2d{\theta}^2 -r^2{sin}^2(\theta)d{\phi}^2$ and $\xi={\xi}^id_i \in T$, let us  introduce ${\xi}_i={\omega}_{ri}{\xi}^r$ with the $4$ {\it formal derivatives} $(d_0=d_t, d_1=d_r,d_2= d_{\theta}, d_3=d_{\phi})$. With speed of light $c=1$ and $A=1-\frac{m}{r}$ where $m$ is a constant, the metric can be written in the diagonal form:  \\
\[ \left(    \begin{array}{cccc}
A & 0  &  0  &  0   \\
0 & -1/A & 0 & 0 \\
0 & 0 & -r^2 & 0 \\
0 & 0 & 0 & -r^2sin^2(\theta)
\end{array}  \right)   \]
Using the notations that can be found in the theory of differential modules, let us consider the Killing equations:   \\
\[ \Omega \equiv {\cal{L}}(\xi)\omega=0 \hspace{1cm} \Leftrightarrow \hspace{1cm} {\Omega}_{ij}\equiv d_i{\xi}_j +d_j{\xi}_i  - 2 {\gamma}^r_{ij}{\xi}_r=0\]
where we have introduced the Christoffel symbols $\gamma$ while  setting $A'={\partial}_rA$ in the differential field $K$ of coefficients ([8], p 87). 
As in the previous Macaulay example and in order to avoid any further confusion between sections and derivatives, we shall use the sectional point of view and rewrite the previous equations in the symbolic form $  L({\xi}_1)\omega = \Omega \in S_2T^* $ where $L$ is the {\it formal Lie derivative}:  \\ 
\[  \left\{ \begin{array}{lcr}
  {\xi}_{3,3} + sin(\theta)cos(\theta){\xi}_2 +rAsin^2(\theta) {\xi}_1 & = & \frac{1}{2}{\Omega}_{33}  \\
  {\xi}_ {2,3}+ {\xi}_ {3,2}  - 2 cot(\theta)   {\xi}_ 3  &  =  & {\Omega}_{23} \\
 {\xi}_ {1,3}+ {\xi}_{3,1}   -  \frac{2}{r}  {\xi}_ 3   &  =  & {\Omega}_{13}   \\
 {\xi}_ {0,3}+ {\xi}_ {3,0}   & = & {\Omega}_{03}    \\
  {\xi}_{2,2} + r A {\xi}_1 & = & \frac{1}{2}{\Omega}_{22}   \\
  {\xi}_{1,2}+ {\xi}_{2,1}  - \frac{2}{r}   {\xi}_ 2  & = & {\Omega}_{12}    \\
 {\xi}_{0,2}+ {\xi}_{2,0}    & = & {\Omega}_{02}   \\
{\xi}_{1,1} + \frac{A'}{2A}{\xi}_1 & = & \frac{1}{2}{\Omega}_{11}   \\
 {\xi}_{0,1} + {\xi}_{1,0} - \frac{A'}{A}{\xi}_0 & = &  {\Omega}_{01}   \\
   {\xi}_{0,0}  - \frac{AA'}{2}{\xi}_1 & = & \frac{1}{2}{\Omega}_{00} 
\end{array} \right.  \]   \\

This system $R_1\subset J_1(T)$ is far from being involutive because it is finite type with second symbol $g_2=0$ defined by the 40 equations $v^k_{ij}=0$ in the initial coordinates. From the symmetry, it is clear that such a system has at least 4 solutions, namely the time translation ${\partial}_t \leftrightarrow {\xi}^0=1 \Leftrightarrow {\xi}_0= A$ and, using cartesian coordinates $(t,x,y,z)$, the 3 space rotations $y{\partial}_z-z{\partial}_y, z{\partial}_x - x{\partial}_z, x{\partial}_y - y{\partial}_x$. \\

These results are brought by the formal Lie derivative of the Weyl tensor because the Ricci tensor vanishes by assumption and we have the splitting $Riemann\simeq Ricci \oplus Weyl$ according to the {\it fundamental diagram II} that we discovered as early as in 1988 ([22]), still not acknowledged though it can be found in ([23],[27],[30-31]). In particular, as the $Ricci$ part is vanishing by assumption, we may identify the $Riemann$ part with the $Weyl$ splitting part as tensors ([31], Th 4.8 and [33]) and it is possible to prove (using a tedious direct computation or computer algebra) that there are only $6$ non-zero components. It is important to notice that this result, bringing a strong condition on the zero jets because of the Lie derivative of the Weyl tensor and thus on the first jets, involves indeed the first derivative of the Weyl tensor because we have a term in $(A'')'$. When 
$\Omega =0$, {\it we obtain after} $2$ {\it prolongations} the additional $5$ new first order PD equations:\\
\[  {\xi}_1=0, {\xi}_{1,2}=0, {\xi}_{1,3}=0, {\xi}_{0,2}=0, {\xi}_{0,3}=0  \]
As we are dealing with sections, ${\xi}_1=0$ {\it does imply} ${\xi}_{1,1}=0$ and ${\xi}_{0,0}=0$ but {\it does not imply} ${\xi}_{1,0}=0$, these later condition being only brought by one additional prolongation and we have the strict inclusions  $R^{(3)}_1\subset R^{(2)}_1 \subset R^{(1)}_1=R_1$. Hence, it remains to determine the dimensions of the subsystems $R'_1=R^{(2)}_1$ and $R"_1=R^{(3)}_1$ with the strict inclusion $R"_1\subset R'_1$, exactly again like in the Macaulay example. Knowing that $dim(R_1)=dim(R_2)=10$, $ dim(R_3) = 5$, $dim(R_4)= 4$, we have thus obtained the $15$ equations defining $R'_1$ with $dim(R'_1)=20-15=5$ and the $16$ equations defining $R"_1$ with $dim(R"_1)=20-16=4$, namely:  \\
\[  \begin{array}{l}
{\xi}_{3,3}+ sin(\theta)cos(\theta) {\xi}_2 =0 \\
{\xi}_{2,3} +{\xi}_{3,2}-2cot(\theta)\, {\xi}_3  =0  \\
{\xi}_{1,3}  =0  \\
{\xi}_{0,3} =0  \\
{\xi}_{2,2} =0 \\
{\xi}_{1,2} =0  \\
{\xi}_{0,2}  =0  \\
{\xi}_{3,1} - \frac{2}{r}\,{\xi}_3 =0  \\
{\xi}_{2,1} - \frac{2}{r} \,{\xi}_2  =0  \\
{\xi}_{1,1} =0  \\
{\xi}_{0,1}  - \frac{A'}{A}\,{\xi}_0=0 \\
{\xi}_{3,0}=0  \\
{\xi}_{2,0}=0  \\
{\xi}_{1,0}=0  \\
{\xi}_{0,0} =0  \\
{\xi}_1=0
\end{array}   \]
Setting now in an intrinsic way ${\xi}_0=A{\xi}^0, {\xi}_1= - \frac{1}{A}{\xi}^1, {\xi}_2= - r^2{\xi}^2$ and in a non-intrinsic way ({\it care}) $ {\xi}_3= - r^2{\xi}^3$, we may even simplify these equations and get a system {\it not depending on} $A$ {\it anymore}: \\
\[  \left\{\begin{array}{l}
{\xi}^3_3+ sin(\theta)cos(\theta) {\xi}^2 =0 \\
{\xi}^2_3 +{\xi}^3_2 - 2cot(\theta)\, {\xi}^3  =0  \\
{\xi}^1_3  =0  \\
{\xi}^0_3 =0  \\
{\xi}^3_1  =0  \\
{\xi}^2_1   =0  \\
{\xi}^1_1 =0  \\
{\xi}^0_1  =0 \\
{\xi}^3_0  =0  \\
{\xi}^2_0=0  \\
{\xi}^1_0=0  \\
{\xi}^0_0 =0  \\
{\xi}^2_2 =0 \\
{\xi}^1_2 =0  \\
{\xi}^0_2  =0  \\
{\xi}^1=0
\end{array}   \right.  \]  \\

It is easy to check that $R^{(3)}_1$, having minimum dimension equal to $4$, is formally integrable, though not involutive as it is finite type, and to exhibit $4$ solutions linearly independent over the constants. Indeed, we must have ${\xi}^0=c$ where $c$ is a constant and we may drop the time variable not appearing elsewhere while using the equation ${\xi}^1=0$. It follows that ${\xi}^2=f(\theta,\phi), {\xi}^3=g(\theta, \phi)$ while $f,g$ are solutions of the first, second and fifth equations of Killing type wih a general solution depending on $3$ constants, a result leading to an elementary probem of $2$-dimensional elasticity left to the reader as an exercise. The system $R^{(3)}_1$ is formally integrable while the system $R^{(2)}_2$ is involutive. Having in mind the PP procedure, it follows that the CC could be of order $2,3$ {\it and even} $4$. Equivalently, we may cut the integration of this system into three systems: \\

\noindent
1) First of all, we have ${\xi}^1=0$ and thus ${\xi}^1_0=0, {\xi}^1_1=0, {\xi}^1_2=0, {\xi}^1_3=0$.  \\
2) Then, we may consider ${\xi}^0_0=0, {\xi}^0_1=0, {\xi}^0_2=0, {\xi}^0_3=0 \Rightarrow  {\xi}^0=c$.  \\
3) Finally, we arrive to the FI system with the same properties as the ones found for Example $2.1$:\\
\[  \left\{\begin{array}{l}
{\xi}^3_3+ sin(\theta)cos(\theta) {\xi}^2 =0 \\
{\xi}^2_3 +{\xi}^3_2 - 2cot(\theta)\, {\xi}^3  =0  \\
{\xi}^3_1  =0  \\
{\xi}^2_1   =0  \\
{\xi}^2_2 =0 \\
\end{array}   \right.  \]  \\
that is with $3$ generating first order CC and $1$ additional second order generating CC.  \\

Proceeding like in the motivating examples, we may introduce the inhomogeneous systems:  \\
\[ \{{\xi}_1=U, {\xi}_{1,2}=V_2, {\xi}_{1,3}=V_3, {\xi}_{0,2}=W_2, {\xi}_{0,3}=W_3\}\in j_2(\Omega) \]
and we finally obtain $16$ PD equations, namely ${\xi}_1=U$ plus the $15$ PD equations:  \\
\[  {\xi}_{0,0}= \frac{1}{2}{\Omega}_{00} + \frac{AA'}{2}U, \,\,\,{\xi}_{0,1} - \frac{A'}{A}{\xi}_0={\Omega}_{01}- \fbox {$d_0U$}\in j_3(\Omega), \,\,\,{\xi}_{0,2}= W_2, \,\,\,{\xi}_{0,3}=W_3 \]
\[  {\xi}_{1,0} = \fbox{$d_0U$}\in j_3(\Omega), \,\,\, {\xi}_{1,1} =\frac{1}{2}{\Omega}_{11}-\frac{A'}{2A}U, \,\,\, {\xi}_{1,2}=V_2, \,\,\, {\xi}_{1,3}=V_3 \]
\[ {\xi}_{2,0}= {\Omega}_{02} -W_2, \,\,\, {\xi}_{2,1} - \frac{2}{r}{\xi}_2= {\Omega}_{12} - V_2, \,\,\, {\xi}_{2,2}= \frac{1}{2}{\Omega}_{22} - rAU \]
\[ {\xi}_{3,0}= {\Omega}_{03} - W_3, \,\,\, {\xi}_{3,1} - \frac{2}{r}{\xi}_3= {\Omega}_{13} -V_3 , \,\,\, {\xi}_{3,2}+{\xi}_{2,3} -
2cot(\theta){\xi}_3={\Omega}_{23}, \]  
\[{\xi}_{3,3}+ sin(\theta)cos(\theta){\xi}_2= \frac{1}{2}{\Omega}_{33} - rAsin^2(\theta) U  \]
As a byproduct, we have $dim(R_3)=5, dim(R_4)=4$ and we obtain $15$ second order CC in $j_2(\Omega)$ along the ker/coker exact sequence:  \\
\[   \begin{array}{lccccccl}
  0  \rightarrow R_3 & \rightarrow & J_3(T) &  \longrightarrow &  J_2(S_2T^*) & \rightarrow &  Q_2  \rightarrow 0 \\
                             &                &  &    \searrow  \hspace{8mm} \nearrow &  &                  &       \\
                             &                 &              &      B_2   &                     &                    &   \\
                             &                 &   &    \nearrow \hspace{8mm}   \searrow     &      &              &      \\
                             &              &     \hspace{8mm}     0        &            &           0  \hspace{10mm}       &                 &
  \end{array}   \]
\[  \begin{array}{rccccccl}
 0  \rightarrow 5 & \rightarrow &  140  &  \rightarrow  &  150 &   \rightarrow &  15 \rightarrow 0   \\ 
  &                &  &    \searrow  \hspace{8mm} \nearrow &  &                  &       \\
                             &                 &              &     135   &                     &                    &  
  \end{array}   \]
 
Then, we have identities to zero like $d_0{\xi}^1-{\xi}^1_0=0$ but we have also {\it surely} the three third order CC like $d_1{\xi}^1- {\xi}^1_1=0, d_2{\xi}^1-{\xi}^1_2, d_3{\xi}^1-{\xi}^1_3$, then {\it  perhaps} the other third order CC $d_2{\xi}^0_3 - d_3{\xi}^0_2=0$ and {\it perhaps} even fourth order CC like $d_0{\xi}^0_1 - d_1{\xi}^0_0=0$ which is containing the leading term $d_{00}U$ after substitution. However, we have the linearization formulas: \\ 
\[  {\rho}_{kl,ij}={\omega}_{kr}{\rho}^r_{l,ij} \,\,\,  \Rightarrow   \,\,\, R_{kl,ij}= {\omega}_{kr}R^r_{l,ij} + {\rho}^r_{l,ij} {\Omega}_{kr}   \] 
\[  {\rho}_{ij}={\rho}^r_{i,rj}\,\,\,  \Rightarrow  \,\,\, {\omega}^{rs}R_{ri,sj} = R_{ij} + {\omega}^{rs}{\rho}^t_{i,rj}{\Omega}_{st}\neq R_{ij} \]
and obtain therefore the formulas:  \\
\[  \begin{array}{rcccl}
    \frac{2r}{3m\,\,sin^2(\theta)}R_{23,30} + \frac{2}{3}{\Omega}_{02} & = \,\,{\xi}_{0,2}\,\, = & - \frac{2r^3A}{3m}R_{01,12}+\frac{1}{3}{\Omega}_{02} &\Leftrightarrow &R_{02}= R^1_{0,12} + R^3_{3,32}=0 \\                                                            
  \frac{2r}{3m}R_{23,02} +\frac{2}{3}{\Omega}_{03} & = \,\, {\xi}_{0,3} \,\, = &  - \frac{2r^3A}{3m}R_{01,13} +\frac{1}{3}{\Omega}_{03} 
   & \Leftrightarrow  &R_{03}= R^1_{0,13} + R^2_{0,23}=0 
\end{array}  \]
with two similar ones for ${\xi}_{1,2}$ and $ {\xi}_{1,3}$ showing the unexpected {\it partition} of the Ricci tensor:  \\
\[   \{ R_{ij}\}=\{ R_{00},R_{11},R_{22},R_{33}\}+ \{ R_{12},R_{13},R_{02},R_{03}\} + \{R_{01},R_{23}\}  \]  
determined by the $15=10 + 5 = 4 + 4 + 2$ second order CC that we have exhibited.  \\
Now, after one prolongation, we get:  \\
\[  {\xi}_{1,00} + (\frac{AA"}{2} - \frac{(A')^2}{4}){\xi}_1+ \frac{1}{2}d_1{\Omega}_{00} - d_0{\Omega}_{01}- \frac{A'}{2A}{\Omega}_{00} + \frac{AA'}{24}{\Omega}_{11} =0     \]
and thus $d_1{\xi}_{0,0} - d_0 {\xi}_{0,1}=0$. Similarly, we have:  \\
\[{\xi}_{1,01} + \frac{A'}{2A}{\xi}_{1,0} - \frac{1}{2}d_0{\Omega}_{11} \equiv d_0 ( d_1{\xi}_1+ \frac{A'}{2A} {\xi}_1- {\Omega}_{11}) = 0  \]
and thus $d_1{\xi}_{1,0} - d_0 {\xi}_{1,1}=0$. It follows that $d_1U + \frac{A'}{2A} - \frac{1}{2}{\Omega}_{11}=0$ is a generating CC of order 
$3$ but $d_{01} U  - \frac{A'}{2A}d_0U - d_0{\Omega}_{11}=0$ is not a generating CC of order $4$. \\
In order to proceed further on, we notice that the generating CC of order $3$ already found can be written as:  \\
\[  d_1U + \frac{A'}{2A}U - \frac{1}{2} {\Omega}_{11}=0, d_2U-V_2=0, d_3U - V_3=0\]
Using crossed derivatives, we get:  \\
\[d_1V_2 +\frac{A'}{2A}d_2U - \frac{1}{2}d_2{\Omega}_{11}=0, d_1V_3 + \frac{A'}{2A}d_3U - \frac{1}{2} d_3{\Omega}_{11}=0, d_2V_3 - d_3V_2=0\]
and thus $d_1{\xi}_{1,2} - d_2{\xi}_{1,1}=0, d_1{\xi}_{1,3} - d_3{\xi}_{1,1}=0, d_2{\xi}_{1,3} - d_3{\xi}_{1,2}=0$.  \\
However, in order to prove that $d_2{\xi}_{0,3} -d_3{\xi}_{0,2}=0$ or equivalently that $ d_2W_3 - d_3W_2=0$, the previous procedure cannot work but we must never forget that $U,V_2,V_3, W_2, W_3$ both belong to $j_2(\Omega)$. Introducing the formal Lie derivative $R=L({\xi}_1)\rho$, we recall that:  \\
\[  W_2= - \frac{2r^3A}{3m}R_{01,12} +\frac{1}{3}{\Omega}_{02}, \hspace{1cm}W_3= - \frac{2r^3A}{3m}R_{01,13} +\frac{1}{3}{\Omega}_{03}  \]
Hence, linearizing the {\it Bianchi identity}:  \\
\[{\nabla}_1{\rho}_{01,23} + {\nabla}_2{\rho}_{01,31} + {\nabla}_3{\rho}_{01,12}=0 \,\,\, \Rightarrow  \,\,\, d_1R_{01,23} + d_2R_{01,31}  + 
d_3R_{01,12} + ... =0  \]
we have proved in ([20]) that the third order CC $d_2W_3-d_3W_2=0$ is not a generating one because it is just a differential consequence of the second order CC $R_{01,23}=0$. \\

Finally, as already noticed, the symbol $g'_1\subset g_1\subset T^*\otimes T$ is not involutive and even $2$-acyclic because otherwise there should only be first order CC for the right members defining the system $R'_1\subset J_1(T)$. As a byproduct, we have, at least on the symbol level, the second 
order CC:  \\
\[  d_{22}{\xi}_{3,3} +d_{33}{\xi}_{2,2} -d_{23}({\xi}_{3,2}+{\xi}_{2,3})=0   \]
and thus:  \\
\[ d_{22}(\frac{1}{2}{\Omega}_{33} -rAsin^2(\theta){\xi}_1 - sin(\theta)cos(\theta){\xi}_2) + d_{33}(\frac{1}{2}{\Omega}_{22} - rA{\xi}_1) - d_{23}({\Omega}_{23}+ 2cot(\theta){\xi}_3)=0  \]
containing surely $d_{22}{\xi}_1, \,\, d_{22}{\xi}_2, \,\, d_{33}{\xi}_1, \,\, d_{23}{\xi}_3$ and thus surely $d_2U,\,\,d_2V_2, d_3V_3$, producing therefore a third order CC that cannot be reduced by means of any Bianchi identity, that is we finally have $15$ generating second order CC and $4$ new generating third order CC, in a manner absolutely similar to that of all the motivating examples of this paper.  \\  

As shown in ([20]), the study of the Killing system for the Kerr metric is even more difficult because the space of solutions is reduced from $4$ already given to the $2$ infinitesimal generators $ \{{\partial}_t, {\partial}_{\phi}\}$ only. Accordingly, we discover that the Schwarzschild and the Kerr metrics do behave quite differently but that {\it there is no hope at all for selecting specific solutions of the Einstein equations in vacuum}. We consider this result as a key challenge when questioning the origin and existence of gravitational waves in general relativity and believe this problem has never been pointed out clearly for the very simple reason that the underlying mathematics are not known by physicists.  \\

\noindent
{\bf EXAMPLE 2.2 REVISITED}: \\

Coming back to the system $R'_2=R^{(2)}_2 \subset R_2$ with a strict inclusion and second members 
$(u,v,w=u_{33}-v_{12}-x^2u_1)$, let us exchange $x^1$ with $x^2$ in order to have an involutive third order symbol $g'_3$ in $\delta$-regular coordinates and consider the system $R'_3$ with now $w=u_{33}-v_{12} - x^1u_2$ in the new coordinates:   \\
\[  \left \{  \begin{array}{lcl}  
 y_{333} - x^1 y_{23}&  =  & v_3  \\
 y_{233} - x^1y_{22} & =  & v_2   \\
 y_{223}  &  =  &  w_3 \\
 y_{222} &  =  &  w_2  \\
 y_{133} - x^1y_{12}& = & v_1  \\
 y_{123}   & = & u_3  \\
 y_{122}  & = & u_2  \\
 y_{112}  & = & u_1  \\ 
  y_{33}-x^1y_2 & = & v \\
  y_{22} &  = &  w   \\
  y_{12}  &  =  &  u
  \end{array} \right. \fbox{ $\begin{array}{ccc}
  1 & 2 & 3  \\
  1 & 2 & \bullet  \\
  1 & 2 & \bullet \\
  1 & 2 & \bullet  \\
  1 & \bullet & \bullet  \\
  1 & \bullet & \bullet  \\
  1 & \bullet & \bullet  \\
  1  &  \bullet & \bullet   \\
  \bullet & \bullet & \bullet  \\
  \bullet & \bullet & \bullet  \\
 \bullet & \bullet & \bullet  
\end{array} $ }   \]
This new system is easily seen to be involutive and we have $3+8+9=20$ first order CC if we consider the second members as just simple notations. Substituting and taking now into account that we have in fact $u_2=d_2u$ formally and so on, all these CC reduce to identities to zero of the form $0=0$ but, using again the original coordinates, $A\equiv w_2-u_1=0, A_2=0, A_3=0 , B\equiv w_{33}- v_{11}-x^2w_1=0$, a system which is {\it not} FI. Accordingly, the generating CC are desribed by $A$ of order $3$ and $B$ of order $4$ with $A_{33}- x^2 A_1 -B_2=0$. I remains to check that this result is coherent with the diagrams of the previous section.  \\

For this, we let the reader compute by hands or with computer algebra the following dimensions $dim(R_2)= 1+3+4=8, \,\, dim(R_{3}= dim(J_3(E)) - dim J_1(F_0)= 12$ because there is no CC of order $1$, $ dim(R_4)) =15 \Rightarrow dim(Q_2)=dim(R_4) - dim(J_4(E)) + dim(J_2(F_0))  = 
15 - 35 + 20 = 0$ because there is no CC of order $2$, $dim(R_5)= 17 \Rightarrow dim(Q_3)= dim(R_5) - dim(J_5(E)) + dim(J_3(F_0)= 
17 - 56 + 40 = 1$ because there is only $1$ CC of order $3$. We have therefore the long sequence:  \\
\[  0 \rightarrow  R_6 \rightarrow J_6(E) \rightarrow J_4(F_0) \rightarrow J_1(Q_2) \rightarrow 0  \]
\[  0 \rightarrow  19 \rightarrow 84 \rightarrow70 \rightarrow 4 \rightarrow 0  \]
and obtain $dim(H(J_4(F_0))= (70 - 4) - (84 - 19)= 66 - 65 = 1$ both with $dim(Q_4)=5$, in a coherent way with the only CC $A$ of order $3$. We let the reader prove that we have similarly $dim(Q_5)=13$ by taking into account the fact that $B_2=A_{33} - x^2A_1$. In order to take into account the existence of a new generating CC of order $4$, we let the reader check that $dim(Q_4)=5 $ and set $F_1=Q_4$ in order to define a fourth order operator ${\cal{D}}_1:F_0 \rightarrow F_1$ by the involutive system:  \\
\[  \left \{ \begin{array}{lcccl}
 B     & \equiv  & u_{3333} + ... & = 0  \\
 A_3 & \equiv  & u_{2333} + ... & = 0 \\
 A_2 & \equiv  & u_{2233} + ... & = 0  \\
 A_1 & \equiv  & u_{1233} + ... & =0  \\
 A  &  \equiv   & u_{233} +  ... & =0
 \end{array} \right. \fbox{ $ \begin{array}{ccc}
 1 & 2 & 3  \\
 1 & 2 & \bullet  \\
 1 & 2 & \bullet  \\
 1 & \bullet & \bullet \\
 \bullet & \bullet & \bullet 
 \end{array} $ }  \]
 Starting anew from this operator, we obtain the first order involutive system:  \\
 \[  \left \{ \begin{array}{lcl}
 d_3A_3 - x^2A_1- B_2 & = & 0 \\
 d_3A_2 - d_2A_3  & =  & 0   \\
 d_3A_1 -d_1A_3 & = & 0  \\
 d_3A - A_3 & = & 0  \\
 d_2A_1 - d_1A_2 & = & 0  \\
 d_2A - A_2 & = & 0  \\
 d_1A - A_1 & = & 0
 \end{array} \right.  \fbox { $ \begin{array}{ccc}
  1 & 2 & 3  \\
  1 & 2 & 3 \\
  1 & 2 & 3 \\
  1 & 2 & 3 \\
 1 & 2 & \bullet  \\
 1 & 2 & \bullet  \\
 1 & \bullet & \bullet 
\end{array} $ }
\Rightarrow  \fbox { $ \begin{array}{ccc}
1 & 2 & 3  \\
1 & 2 & 3 \\
1 & 2 & 3 \\
1 & 2 & \bullet 
 \end{array} $ }
 \Rightarrow  \fbox  { $ \begin{array}{ccc}
 1 & 2 & 3
 \end{array}  $ }  \]
 where each Janet tabular is induced from the preceding one till the end of the procedure as in ([12], p 153,154  for details). We also notice that this system 
brings {\it automatically} the Spencer operator. \\
We obtain therefore the following differential sequence:  \\
\[  0 \rightarrow \Theta \rightarrow E \underset 2 {\stackrel{{\cal{D}}}{\longrightarrow}}F_0 \underset 4 {\stackrel{{\cal{D}}_1}{\longrightarrow}} F_1 \underset 1 {\stackrel{{\cal{D}}_2}{\longrightarrow}} F_2 \underset 1 {\stackrel{{\cal{D}}_3}{\longrightarrow}} F_3 
\underset 1 {\stackrel{{\cal{D}}_4}{\longrightarrow}} F_4 \rightarrow 0  \]
\[  0 \rightarrow \Theta \rightarrow 1 \underset 2 {\longrightarrow}2\underset 4 {\longrightarrow} 5 \underset 1 {\longrightarrow} 7
\underset 1 {\longrightarrow} 4 \underset 1 {\longrightarrow} 1 \rightarrow 0  \]
which is formally exact on the jet level, even if it is not strictly exact because the first operator is not FI, and we check that $1-2+5-7+4-1=0$. We notice that the part between $F_0$ and $F_4$ is typically a Janet sequence for ${\cal{D}}_1$. \\
It follows that we have the following long exact sequence on the level of jets, $\forall r\geq -5$:  \\
\[   0 \rightarrow R_{r+9}\rightarrow J_{r+9}(E) \rightarrow J_{r+7}(F_0) \rightarrow J_{r+3}(F_1) \rightarrow J_{r+2}(F_2) \rightarrow J_{r+1}(F_3) \rightarrow J_r(F_4) \rightarrow 0  \]
a result leading to:   \\
\[  \begin{array}{lcl}
 dim(R_{r+9})  & = &  1(r+10)(r+11)(r+12)/6 - 2(r+8)(r+9)(r+10)/6 \\
                        &    & + 5 (r+4)(r+5)(r+6)/6 - 7(r+3)(r+4)(r+5)/6  \\
                        &   & + 4(r+2)(r+3)(r+4)/6 -1(r+1)(r+2)(r+3)/6     \\
                        & = &  2r + 25  
                        \end{array}  \]
and thus to $dim(R_{r+4})=2r + 15, \forall r\geq 0$, a result not evident to grasp at first sight because it comes from the lack of formal integrability of $R_2$ and the strict inclusion $R^{(2)}_2 \subset R_2$. \\   \\

\newpage 
\noindent
{\bf EXAMPLE 2.3 REVISITED}:  \\

Coming back to the systems $R_2$ with second members $(u,v)$ and $R'_2=R^{(2)}_2 \subset R_2$ with a strict inclusion and second members 
$(u,v,2w=u_{33}-v_{22}-x^2u_1)$, let us exchange $x^1$ with $x^2$ in order to have an involutive third order symbol $g'_3$ in $\delta$-regular coordinates but the system $R'_3$, with now $w=u_{33}-v_{12} - x^1u_2$ in the new coordinates, is {\it not} FI. Hence, we must start anew with the system $R"_2=R^{(4)}_2 \subset R'_2$ with a strict inclusion, described by the $4$ PD equations:  \\
\[  \left\{  \begin{array}{lcl}
  y_{33}- x^2 y_1  & =  &  v  \\
 y_{22}  & = & u  \\
 y_{12}   &  =  &  w  \\
 y_{11}  &   =  &  \omega
 \end{array}   \right.  \]
where $2\omega = u_{3333}- v_{2233} - 2 x^2 u_{133} + x^2 v_{122} - 2 v_{12} + (x^2)^2u_{11}$. Using one prolongation, we get the third order PD equations:   \\
\[  \left \{ \begin{array}{lcl}
 y_{333} - x^2 y_{13}& = & v_3 \\
 y_{233} - y_1 &  =  &  v_2 + x^2 w   \\
y_{223}& = & u_3  \\
 y_{222} & = & u_2  \\
y_{133}   & = & v_1 + x^2 \omega \\
y_{123} & = & w_3 \\
y_{122} & = & u_1 = w_2 \Rightarrow A   \\
y_{113} &  =  &  {\omega}_3  \\
y_{112}  &  =  &  w_1={\omega}_2 \Rightarrow B  \\
y_{111} & =  & {\omega}_1
 \end{array} \right.  \fbox { $ \begin{array}{ccc}
  1 & 2 & 3  \\
  1 & 2 & \bullet \\
  1 & 2 & \bullet \\
  1 & 2 & \bullet \\
 1 & \bullet & \bullet  \\
 1 & \bullet & \bullet  \\
 1 & \bullet & \bullet  \\
 1 & \bullet &  \bullet  \\
 1 & \bullet & \bullet  \\
 1 & \bullet & \bullet 
\end{array} $ }   \]
and we discover that the symbol $g"_2$ is finite type because $g"_3=0$. As we had to use one prolongation in order to get a $2$-acyclic symbol, we obtain {\it sixth} order CC $(A,B,C) \in j_2(u,v,w,\omega)$. We refer the reader to ([13], p 83) or ([14], p 315) for more details on this delicate result.

Using the notations of the last section, we now provide the systems $B_2, B_3,B_4, B_5, B_6$ together and we notice the following striking results:  \\
\[  B_2=J_2(F_0), B_3 \subset {\rho}_1(B_2)\subset J_3(F_0), B_4 ={\rho}_1(B_3)\subset J_4(F_0), \]
\[   B_5={\rho}_1(B_4) \subset J_5(F_0), B_6 \subset {\rho}_1(B_5) \subset J_6(F_0)  \]
reaching therefore the following involutive system of order $6$ where we did not quote $B=0$ because we already proved that $B= A_{33} - x^2 A_1$:\\

\[  \left \{ \begin{array}{lclcl}
 C     & \equiv  & u_{333333} + ... & = 0  \\
 A_{333} & \equiv  & u_{233333} + ... & = 0 \\
 A_{233} & \equiv  & u_{223333} + ... & = 0  \\
 A_{223} & \equiv  & u_{222333}  + ... & =0  \\
 A_{222} & \equiv  & u_{222233}  + ...  & = 0  \\
 A_{133} & \equiv  & u_{123333} + ... & = 0  \\
 A_{123} & \equiv  & u_{122333} + ...  &  = 0  \\
 A_{122} & \equiv  & u_{122233} + ...  &  =  0  \\
 A_{113} & \equiv  & u_{112333} + ...  & =  0  \\
 A_{112} & \equiv  & u_{112233}  + ... &  = 0  \\
 A_{111} & \equiv  & u_{111233} + ...  & =0   \\                                                                                                                                                                                                                                                                                                                                                                                                                                                                                                                                                                                                                                                                                                                                                         
 A_{33}   & \equiv  & u_{23333} + ... & = 0 \\
 A_{23}   & \equiv  & u_{22333} + ...  & = 0  \\
 A_{22}   & \equiv  & u_{22233} + ...  & = 0  \\
 A_{13}   & \equiv  & u_{12333} + ...  & = 0  \\
 A_{12}   & \equiv & u_{12233} + ... & = 0  \\
 A_{11}   & \equiv  & u_{11233} +  ... & = 0  \\
 A_3        & \equiv  & u_{2333} + ...  & =  0  \\
 A_2        & \equiv  & u_{2233} + ...  & = 0  \\
 A_1        & \equiv  & u_{1233} + ...  & =  0  \\
 A           &  \equiv & u_{233} +  ... & =0
 \end{array} \right. \fbox{ $ \begin{array}{ccc}
 1 & 2 & 3  \\
 1 & 2 & \bullet  \\
 1 & 2 & \bullet  \\
 1  &  2  &  \bullet  \\
 1  &  2  &  \bullet  \\
 1 & \bullet & \bullet \\
 1  & \bullet & \bullet \\
 1 & \bullet  & \bullet  \\
 1  &  \bullet  &  \bullet  \\
 1  &  \bullet  & \bullet  \\
 1  &  \bullet  &  \bullet  \\
 \bullet & \bullet & \bullet \\
  \bullet & \bullet & \bullet \\
  \bullet & \bullet & \bullet \\
  \bullet & \bullet & \bullet \\
 \bullet & \bullet & \bullet \\
  \bullet & \bullet & \bullet \\
  \bullet & \bullet & \bullet \\
  \bullet & \bullet & \bullet \\
 \bullet & \bullet & \bullet \\
  \bullet & \bullet & \bullet 
 \end{array} $ }  \]
 Starting anew from this operator providing $46$ CC, we obtain the first order involutive system:  \\
 \[  \left \{ \begin{array}{lcl}
 d_3A_{333} - 2 x^2A_{133} + (x^2)^2A_{11} - C_2 & = & 0 \\
 d_3A_{233}- d_2A_{333} & =  & 0   \\
 d_3A_{223} - d_2A_{233} &  =  &  0  \\
 d_3A_{222} - d_2A_{223} &  =  &  0  \\
 .......................................&.....&.....\\
 d_3A_1 -d_1A_3 & = & 0  \\
 d_3A - A_3 & = & 0  \\
 d_2A_{133} - d_1A_{233} & = &  0  \\
 ...............   &..&... \\
 d_2A_1 - d_1A_2 & = & 0  \\
 d_2A - A_2  &  =  &  0     \\
 d_1A_{33} - A_{133} & = & 0  \\
 ...........................&....& ...  \\
 d_1A_1 - A_{11} & = & 0  \\
 d_1A - A_1 & = & 0
 \end{array} \right.  \fbox { $ \begin{array}{ccc}
  1 & 2 & 3  \\
  1 & 2 & 3 \\
  1 & 2 & 3 \\
  1 & 2 & 3 \\
  ...  &...  & ...   \\
  1  &  2  &   3  \\
  1 & 2 & 3 \\
   1 & 2 &\bullet  \\
   ...  & ...   &...  \\
    1  &  2  &  \bullet \\
    1  &  2  & \bullet  \\
 1 & \bullet & \bullet  \\
  ... &...  &  ...  \\
 1 &  \bullet & \bullet  \\
 1 & \bullet & \bullet 
\end{array} $ }\]
with $20$ equations of class $3$, $16$ equations of class $2$ and $10$ equations of class $1$. There are $36$ CC  providing an involutive system with $26$ equations of class $3$, $10$ equations of classs $2$ but no equation of class $1$. We get a final system of $10$ CC of class $3$ without any CC. Like in the preceding application, we have thus obtained the following formally exact sequence:  \\
\[  0 \rightarrow \Theta \rightarrow E \underset 2 {\stackrel{{\cal{D}}}{\longrightarrow}}F_0 \underset 6 {\stackrel{{\cal{D}}_1}{\longrightarrow}} F_1 \underset 1 {\stackrel{{\cal{D}}_2}{\longrightarrow}} F_2 \underset 1 {\stackrel{{\cal{D}}_3}{\longrightarrow}} F_3 
\underset 1 {\stackrel{{\cal{D}}_4}{\longrightarrow}} F_4 \rightarrow 0  \]
\[  0 \rightarrow \Theta \rightarrow 1 \underset 2 {\longrightarrow}2\underset 6 {\longrightarrow} 21 \underset 1 {\longrightarrow} 46
\underset 1 {\longrightarrow} 36 \underset 1 {\longrightarrow} 10 \rightarrow 0  \]
with $1 - 2 + 21 - 46 + 36 - 10 = 0$, a part of it being a Janet sequence as before. Similarly, we get:  \\
\[  \begin{array}{lcl}
 dim(R_{r+11})  & = &  1(r+12)(r+13)(r+14)/6 - 2(r+10)(r+11)(r+12)/6 \\
                        &  &    21(r+4)(r+5)(r+6)/6 - 46(r+3)(r+4)(r+5)/6   \\
                        &   & + 36(r+2)(r+5)(r+6)/6 -10(r+1)(r+2)(r+3)/6     \\
                        & = &   18
                        \end{array}  \]
or even $dim(R_{r+6})=18, \forall r\geq 0$ as a striking result indeed that can be checked directly through the exact sequences:  \\
\[    0 \rightarrow R_8 \rightarrow J_8(E) \rightarrow  J_6(F_0) \rightarrow F_1 \rightarrow 0  \]
\[    0 \rightarrow 18 \rightarrow 165 \rightarrow  168\rightarrow 21 \rightarrow 0  \]
\[  .................................................................................\]
\[    0 \rightarrow R_6 \rightarrow J_6(E) \rightarrow  J_4(F_0) \rightarrow Q_4 \rightarrow 0  \]
\[    0 \rightarrow 18 \rightarrow 84 \rightarrow 70  \rightarrow 4 \rightarrow 0  \]
\[   0  \rightarrow R_5  \rightarrow  J_5(E) \rightarrow  J_3(F_0) \rightarrow Q_3  \rightarrow 0  \]
\[  0  \rightarrow 17  \rightarrow 56 \rightarrow 40 \rightarrow  1  \rightarrow 0  \] 
and comes from the fact that ${dim}_K(R)=6 < \infty$ or, equivalently, ${dim}_K(M)=6 < \infty$.      \\
The reader not familiar with the formal theory of differential systems or modules may be surprised by the fact the two dimensions just found do not coincide at all because $6 < 18$. However, we have indeed $dim (R"_3)=dim(R"_2)= dim(R^{(4)}_2)= 6$ and the exact sequence:  \\
\[   0  \rightarrow g_4  \rightarrow R_4  \rightarrow R_3  \rightarrow R_3/R^{(1)}_3  \rightarrow  0  \]
\[  0  \rightarrow 4  \rightarrow 15  \rightarrow 12 \rightarrow 1 \rightarrow  0 \hspace{1cm}\Rightarrow \hspace{1cm} dim(R^{(1)}_3)= 12 - 1 = 11 \]
showing that $dim(R^{(1)}_3) < dim(R_3) < dim(R_4)$ with $ 11 < 12 < 15$.However, we have the general Theorem 2.A.7 in ([20]) providing the useful {\it prolongation/ projection} (PP) procedure, namely that we have ${\rho}_r(R^{(1)}_q)= R^{(1)}_{q+r}, \forall r\geq 0$ whenever the symbol $g_q$ of $R_q$ is $2$-acyclic. In the present case, we have indeed ${\rho}_r(R^{(1)}_{3})= R^{(1)}_{r+3}, \forall r \geq 0$ because $g_3$ is known to be involutive, and the final system $R^{(4)}_3$ is involutive with zero symbol, providing $R^{(4)}_2$ which is only FI but with dimension $6$. This situation is quite tricky indeed because prolongations are filling up successively the PD equations of order $2$, then $3$ and so on, adding therefore:  \\
 \[   \{y_{12}=0\},\{ y_{123}=0, y_{122}=0, y_{112}=0\},\{ y_{11}=0\}, \{y_{113}=0, y_{112}=0, y_{113}=0\}, ...\]     \\

\vspace{15mm}
\noindent
{\bf 5) CONCLUSION}  \\

When a differential operator ${\cal{D}}$ of order $q$ is given, the problem of finding its {\it compatibility conditions} (CC) is to look for a new 
operator ${\cal{D}}_1$ of a certain order $s$ such that ${\cal{D}}_1\eta=0$ {\it must} be satisfied in order to be able to solve the inhomogeneous system $ {\cal{D}}\xi=\eta$. This is an old problem first solved as a footnote by M. Janet in 1920 ([6],[12],[13],[14]) and finally studied by D.C. Spencer in 1970 ([4],[5],[25]). The main idea is to construct a {\it finite length differential sequence} by repeating this procedure anew with ${\cal{D}}_1$ and so on till one eventually ends with ${\cal{D}}_n$ according to Janet when $n$ is the number of independent variables. It soon became clear that constructing ${\cal{D}}_1$ is largely depending on intrinsic properties of ${\cal{D}}$. \\

\noindent
$\bullet $ If ${\cal{D}}$ is {\it involutive}, then ${\cal{D}}_1, ... , {\cal{D}}_n$ are first order involutive operators in the corresponding {\it Janet sequence} that can be constructed "{\it step by step} " as above but also "{\it as a whole}" like in the Poincar\'{e} sequence for the exterior derivative. \\

\noindent
$\bullet$ If ${\cal{D}}$ is only {\it formally integrable} (FI), that is {\it all} the equations of order $q+r$ of the corresponding homogeneous system can be obtained by {\it only} $r$ prolongations, then the order of ${\cal{D}}_1$ is $s+1$ when $s$ is the smallest integer such that the symbol of order 
$q+s$ becomes $2$-acyclic. Such a result is still not acknowledged today by physicists even though it is essential for studying the conformal Killing system of space-time in general relativity.  \\

\noindent
$\bullet$ If ${\cal{D}}$ is not even FI, not only the construction of ${\cal{D}}_1$ may become very difficult but also a {\it strange phenomenon} may appear, namely one can start to find CC of order $s_1$, then no new CC other than the ones generated by {\it these} CC up to order $s_2>s_1$ when {\it suddenly} new generating CC may appear, generating all the CC up to order $s_3>s_2$ and so on till the procedure ends.  \\

This delicate question has been recently raised by physicists in order to study the Killing systems for various useful metrics solutions of the Einstein equations in vacuum (Minkowski gives $s=2$ while Schwarzschild gives $s_1=2,s_2=3$). Needless to say that computer algebra is of quite a poor help in this case because the dimensions of the jet spaces and the size of the matrices involved may increase drastically ([21]).  \\

The aim of this paper has been first to provide illustrating examples of the above situations and one of them with $s_1=3, s_2=6$ seems to be the only one known in the literature today. In addition, we have solved the (general) generating problem by using new differential homological algebraic methods, with the hope that computer algebra will soon become of some help in a near future.  \\

\vspace{4cm}

\noindent
{\bf 6) REFERENCES}  \\

\noindent
[1] Aksteiner, S., B\"{a}ckdahl, T.: All Local Gauge Invariants for Perturbations of the Kerr Spacetime, 2018, arxiv:1803.05341.   \\
\noindent
[2] Andersson, L.: Spin Geometry and Conservation Laws in the Kerr Spacetime, 2015, arxiv:1504.02069.  \\
\noindent
[3] Cartan,E., Einstein,A.: Letters on Absolute Parallelism, Princeton University Press, 1979.  \\
\noindent
[4] Goldschmidt, H.: Prolongations of Linear Partial Differential Equations: I A Conjecture of Elie Cartan, Ann. Scient. Ec. Norm. Sup., 4, 1 (1968) 417-444.  \\
\noindent
[5] Goldschmidt, H.: Prolongations of Linear Partial Differential Equations: II Inhomogeneous Equations, Ann. Scient. Ec. Norm. Sup., 4, 1 (1968) 617-625.  \\
\noindent
[6] Janet, M.: Sur les Syst\`{e}mes aux D\'{e}riv\'{e}es Partielles, Journal de Math., 8(3) (1920) 65-151.  \\
\noindent
[7] Kashiwara, M.: Algebraic Study of Systems of Partial Differential Equations, M\'emoires de la Soci\'{e}t\'{e} Math\'ematique de France 63, 1995, 
(Transl. from Japanese of his 1970 Master's Thesis).\\
\noindent
[8] Khavkine, I.: The Calabi Complex and Killing Sheaf Cohomology, J. Geom. Phys., 113 (2017) 131-169. \\
arXiv:1409.7212  \\
\noindent
[9] Khavkine, I.: Compatibility Complexes of Overdetermined PDEs of Finite Type, with Applications to to the Killing Equation, 2018. \\
arxiv:1805.03751  \\
\noindent
[10] Macaulay, F.S.: The Algebraic Theory of Modular Systems, Cambridge, 1916.  \\
\noindent  
[11] Northcott, D.G.: An Introduction to Homological Algebra, Cambridge university Press, 1966.  \\
\noindent
[12] Pommaret, J.-F.: Systems of Partial Differential Equations and Lie Pseudogroups, Gordon and Breach, New York, 1978; Russian translation: MIR, Moscow, 1983.\\
\noindent
[13] Pommaret, J.-F.: Partial Differential Equations and Group Theory, Kluwer, 1994.\\
http://dx.doi.org/10.1007/978-94-017-2539-2    \\
\noindent
[14] Pommaret, J.-F.: Partial Differential Control Theory, Kluwer, Dordrecht, 2001 (1000 pp).\\
\noindent  
[15] Pommaret, J.-F.: Relative Parametrization of Linear Multidimensional Systems, Multidim. Syst. Sign. Process. (MSSP), Springer, 26 
(2013) 405-437. \\
http://dx.doi.org/10.1007/s11045-013-0265-0  \\
\noindent
[16] Pommaret, J.-F.: The Mathematical Foundations of General Relativity Revisited, Journal of Modern Physics, 4 (2013) 223-239.\\
http://dx.doi.org/10.4236/jmp.2013.48A022   \\
\noindent
[17] Pommaret, J.-F.: Deformation Theory of Algebraic and Geometric Structures, Lambert Academic Publisher, (LAP), 
Saarbrucken, Germany, 2016.  \\
http://arxiv.org/abs/1207.1964  \\
\noindent
[18] Pommaret, J.-F.: Why Gravitational Waves Cannot Exist, Journal of Modern Physics, 8,13 (2017) 2122-2158.  \\
http://dx.doi.org/10.4236/jmp.2017.813130   \\
\noindent
[19] Pommaret, J.-F.: New Mathematical Methods for Physics, NOVA Science Publisher, New York, 2018.  \\
\noindent  
[20] Pommaret, J.-F.: Minkowski, Schwarzschild and Kerr Metrics Revisited, Journal of Modern Physics, 9 (2018) 1970-2007.  \\
https://doi.org/10.4236/jmp.2018.910125   \\
arXiv:1805.11958v2   \\
\noindent
[21] Quadrat, A., Robertz, R.: A Constructive Study of the Module Structure of Rings of Partial Differential Operators, Acta Applicandae Mathematicae, 133 (2014) 187-234. \\
http://hal-supelec.archives-ouvertes.fr/hal-00925533   \\
\noindent
[22] Rotman, J.J.: An Introduction to Homological Algebra, Pure and Applied Mathematics, Academic Press, 1979.  \\
\noindent
[23] Schneiders, J.-P.: An Introduction to D-Modules, Bull. Soc. Roy. Sci. Li\'{e}ge, 63 (1994) 223-295.  \\
\noindent
[24] Shah, A.G., Whitting, B.F., Aksteiner, S., Andersson, L., B\"{a}ckdahl, T.: Gauge Invariant perturbation of Schwarzschild Spacetime, 2016, 
arxiv:1611.08291  \\
\noindent
[25] Spencer, D.C.: Overdetermined Systems of Partial Differential Equations, Bull. Amer. Math. Soc., 75 (1965) 1-114.\\
\noindent

\end{document}